\documentclass{amsart}
\font\emailfont=cmtt10

\usepackage{amsmath,amssymb,amsthm,amsfonts,amscd,flafter,epsf, epsfig,
graphicx}
\usepackage{amssymb,latexsym}

\ifx\mathbb\let\mathbb\Bbb\fi
\newcommand\Cc{\mathbb{C}}
\newcommand\Rr{\mathbb{R}}
\newcommand\Qq{\mathbb{Q}}
\newcommand\Zz{\mathbb{Z}}
\newcommand\Ff{\mathbb{F}}

\newcommand{\conj}{{\rm conj}}
\newcommand{\lk}{{\rm lk}}
\newcommand{\tb}{{\rm tb}}
\newcommand{\sing}{{\rm sing}}
\newcommand{\sitare}{{\rm star}}
\newcommand{\odd}{{\rm odd}}
\newcommand{\even}{{\rm even}}
\newcommand{\pr}{{\rm pr}}

\theoremstyle{plain}
\newtheorem{theorem}{Theorem}
\newtheorem{lemma}[theorem]{Lemma}
\newtheorem{corollary}[theorem]{Corollary}
\newtheorem{proposition}[theorem]{Proposition}
\newtheorem{remark}[theorem]{Remark}

\begin{document}

\title[{Generalised Thurston-Bennequin invariants for real
algebraic surface singularities}] {Generalised Thurston-Bennequin
invariants for real algebraic surface singularities}

\author[FER\.{I}T \"{O}zt\"urk]{FER\.{I}T \"{O}zt\"urk}
\address{Department of
Mathematics, Bo\u{g}azi\c{c}i University, \.{I}stanbul, Turkey
\newline \indent{\emailfont{ferit.ozturk@boun.edu.tr}}}




\begin{abstract}
A generalised Thurston-Bennequin invariant for a Q-singularity of
a real algebraic variety is defined as a {\it linking form} on the
homologies of the real link of the singularity. The main goal of
this paper is to present a method to calculate the linking form in
terms of the very good resolution graph of a real normal unibranch
singularity of a real algebraic surface. For such singularities,
the value of the linking form is the Thurston-Bennequin number of
the real link of the singularity. As a special case of unibranch
surface singularities, the behaviour of the linking form is
investigated on the Brieskorn double points $x^m+y^n\pm z^2=0$.
\end{abstract}

\maketitle

\section{Introduction}
\label{intro}
The link of a real algebraic surface singularity is a contact
3-manifold with the canonical contact structure the set of complex
tangencies. Furthermore, the real link in this contact manifold is
a Legendrian link. The {\it linking form} is defined as a rational
valued form on the first homology of the real link whenever the
rational linking number is well-defined in the link of the
singularity, i.e. when the link is a $\Qq$-homology sphere (see
Section \ref{sec:2} or \cite{Fin1} for the definition of the
linking form). For such a link, the linking form, when evaluated
on a pair of classes of knots of the real link, gives the rational
linking number of the knots with contact framing. In the case when
the real link is connected, i.e. the singularity is unibranch, the
outcome of the linking form is by definition the
Thurston-Bennequin number.

Therefore, one can consider the linking form
as a {\it generalised} Thurston-Bennequin invariant. It is
generalised in the sense that (i) it is defined for a {\it link} in a contact
manifold, (ii) the knots of the link are not necessarily
integrally null-homologous but they bound rationally, (iii) it is rational valued.

A method to compute this linking form has already been presented for
surface singularities $f(x,y)\pm z^2=0$ ($f(x,y)$ reduced) using real deformations of
the singularity \cite{Fin2}.

The main result of this paper is Theorem \ref{calcul}, Section
\ref{sec:6}, which presents a method to calculate the generalised
Thurston-Bennequin number of the connected real link of a normal
real algebraic surface Q-singularity using a very good resolution
graph of the singularity with the condition that the singularity
is numerically Gorenstein. In particular, Theorem \ref{calcul} is
valid for all the unibranch hypersurface singularities in $\Cc^3$.
Furthermore, we give an example which exhibits how the method
works algorithmically for singularities with more than one real
branch.

In the proof of Theorem \ref{calcul}, we make use of the fact that
the generalised Thurston-Bennequin number of the connected real
link of a surface singularity gives by definition the
self-intersection of the oriented real part of the surface in the
complex surface relative to the link of the singularity. This is
not valid anymore if the real part is non-oriented. Theorem
\ref{calcul} determines how much the generalised
Thurston-Bennequin number differs from the self-intersection of
the real part. This difference is a non-integer rational number in
general.


In the second part of the article, as a special case of
hypersurface singularities, for which Theorem \ref{calcul} is
valid, we deduce some properties for the Thurston-Bennequin
numbers of the Brieskorn double points $x^m+y^n\pm z^2=0$ with
$m$, $n$ relatively prime. We first express the conditions under
which the Thurston-Bennequin numbers become integer (see Corollary
\ref{q-int} and Corollary \ref{q+int}, Section \ref{sec:8}). Then
we investigate the behaviour of the Thurston-Bennequin numbers
with changing $m$ and $n$ (see Theorems \ref{period},
\ref{simetri} and \ref{simetreven}, Section \ref{sec:8}). These
theorems show that one can calculate the Thurston-Bennequin
numbers for large $m$ and $n$ in terms of the Thurston-Bennequin
numbers corresponding to sufficiently small $m$ and $n$

The main approach in the proofs is to use Theorem \ref{calcul}. In
fact, Corollary \ref{q-int} is a direct consequence of Theorem
\ref{calcul}. We prove Corollary \ref{q+int} and Theorem
\ref{period} by observing the change in the very good resolution
graph of the singularity when $m$ is held fixed and $n$ is varied
(see Section \ref{sec:4}) and using the results of Section
\ref{sec:5}.

The proofs of Theorems \ref{simetri} and \ref{simetreven}, Section
\ref{sec:8}, do not involve resolution graphs. Rather, we
construct new real algebraic surfaces as branched double coverings
of appropriate weighted homogeneous surfaces; these new algebraic
surfaces contain at the same time both singularities involved in
each theorem.

{\it Notation:} In this article, the manifold $\overline{X}$ is
usually the quotient of the manifold $X$ by complex conjugation
except in the case $\overline{\Cc P^2}$ which is $\Cc P^2$ with
reverse orientation.

\noindent {\bf Acknowledgments.} I would like to thank Sergey
Finashin (METU, Ankara) for having introduced me the problem, for
his invaluable guidance and his inspiration and to Y{\i}ld{\i}ray
Ozan (METU, Ankara) for many stimulating discussions. In addition,
I would like to thank CMAT, Ecole Polytechnique, Palaiseau where I
completed this work as a postdoc.

\section{Preliminary notions and definitions}
\label{sec:0}
\subsection{The real branches of a real surface Q-singularity
and the linking form} \label{sec:1}
Let $\Cc X$ be the complex point set of a  real algebraic surface
with complex conjugation $\conj:\Cc X\rightarrow \Cc X$; let $\Rr
X$ denote the set of the real points of $\Cc X$ and $\overline{X}$
denote the quotient $\Cc X/\conj$. Let $x$ be an isolated
singularity of $\Cc X$ and $\Cc U$ be a small compact cone-like
neighbourhood of $x$ \cite{Mil}.

It is known that the boundary $\Cc M=\partial\,\Cc U$ is a
3-manifold with a canonical contact structure determined by the
distribution of complex lines that exist uniquely in the tangent
space over each point of $\Cc M$ \cite{Var}. Furthermore, since
the real link $\Rr M$ is tangent to the unique complex line at
each point, $\Rr M$ is a Legendrian link in $\Cc M$. Although $\Cc
M$ is not a complex manifold, we use the notation with prefix
$\Cc$ to emphasize that $\Cc M$ is fixed under $\conj$ of $\Cc X$.
With the same idea in mind, $\Rr M$ denotes the subset of $\Cc M$
fixed {\it pointwise} by $\conj$. Finally $\overline{M}$ will
denote $\partial \overline{U}$.

A singular point $x\in\Cc X$ is called a {\em topologically
rational} singularity (or a {\em Q-singularity}) if its link $\Cc
M$ is a $\Qq$-homology sphere. Point $x$ is called a {\em
$\overline{Q}$-singularity} if $\overline{M}$ is a $\Qq$-homology
sphere. Let $x\in\Cc X$ be a Q-singularity and let $\sigma$ denote
the number of connected components of the real link $\Rr M$. Each
connected component of $\Rr M$ is called a {\it real branch} of
$\Cc X$ at $x$. If $\sigma=1$, i.e. if $\Rr M$ is a knot, then we
will say that $x$ is a {\it real unibranch singularity}. For
example, the Brieskorn singularity $\{x^m+y^n+z^d=0\}$,
$\gcd(m,n,d)=1$ is an isolated Q-singularity which is a real
unibranch singularity.

Let $\gamma_1,\ldots,\gamma_{\sigma}$ be the connected components
of $\Rr M$ and $\xi_i$ be a {\it Thurston-Bennequin framing} of
$\gamma_i$ in $\Cc M$, i.e. a choice of trivialisation of a
tubular neighbourhood of $\gamma_i$ in $\Cc M$. Let us denote by
$\gamma_i^{\xi_i}$ the knot obtained by a small shift of
$\gamma_i$ in the direction of $\xi_i$. Consider the bilinear form
\[\tb_x:H_1(\Rr M;\Qq)\times H_1(\Rr M;\Qq)\rightarrow\Qq\]
defined on the generators by:
\[\tb_x([\gamma_i],[\gamma_j])=\lk(\gamma_i,\gamma_j^{\xi_j}) \mbox{ in } \Cc M \]
and extended linearly over $H_1(\Rr M;\Qq)$ (see \cite{Fin1} or
\cite{Fin2}). Here lk is the rational linking number in $\Cc M$.
Note that if $i=j$, then
$\tb_x([\gamma_i],[\gamma_i])=\tb(\gamma_i)$, the rational
Thurston-Bennequin number of the Legendrian knot \cite{Ben}. The
form $\tb_x$ is called the {\em linking form}.

It is essential to note the following: Let $\sigma_i$ denote the
component of closure of $\Rr U - \{x\}$ containing $\gamma_i$; the
real surface $\sigma_i$ is a cone over $\gamma_i$ topologically.
Consider a tangent vector field $\nu_i$ on $\sigma_i$ directed
outward along $\gamma_i$. Then $\sqrt{-1}\cdot\nu_i$ coincides
with $\xi_i$, the Thurston-Bennequin framing, proving that in case
$\sigma_i$ and $\sigma_j$ are oriented, the form $\tb_x$ evaluated
on $[\gamma_i]$ and $[\gamma_j]$ is nothing but the rational
intersection of $\sigma_i$ and $\sigma_j$ in $\Cc U$ relative to
the boundary, i.e. :
\begin{equation} \tb_x([\gamma_i],[\gamma_j])=\langle
\sigma_i,\sigma_j\rangle_{(\Cc U,\Cc M)}.
\label{selfint}\end{equation}

Similarly as above, for a $\overline{Q}$-singularity,
one can define a bilinear form $\overline{\tb}_x$
on $H_1(\Rr M;\Qq)$ by letting
\[\overline{\tb}_x([\gamma_i],[\gamma_j])=\lk(q(\gamma_i),q(\gamma_j^{\xi_j}))
\mbox{ in } \overline{M}=q(\Cc M)\] where $q:\Cc
X\rightarrow\overline{X}$ is the quotient map. For a
$\overline{Q}$-singularity, $\tb_x$ and $\overline{\tb}_x$ satisfy
\cite{Fin2}:
\begin{equation}
\overline{\tb}_x=2\,\tb_x.
\label{tbbar2rb}
\end{equation}

The linking form $\tb_x$ corresponding to
a singularity with $\sigma$ real branches
is represented by
a $\sigma\times \sigma$ matrix. We are going to denote by $[\tb_x]$ the
$\sigma\times \sigma$ matrix in the basis $\gamma_1,\ldots,\gamma_{\sigma}$.
We will not mention $\gamma_i$s whenever there is no confusion and we will
also allow permuted indexing.
The matrix
$[\tb_x]$ of a real unibranch singularity $x$ is
a $1\times 1$ matrix and in that case,
$\tb_x$ will be called
{\it the Thurston-Bennequin number} of the real unibranch singularity $x$
as well.


An algorithmic method was introduced in \cite{Fin2}, Theorem B, to
compute the Thurston-Bennequin numbers of the real unibranch
singularities $f(x,y)\pm z^2=0$ where $f(x,y)$ is a real
polynomial with an isolated real singularity. This method makes
use of the diagrams of Gusein-Zade \cite{Gus} and A'Campo
\cite{ACa} which are obtained by small real generic deformations
of the singularity $f(x,y)=0$ in $\Cc^2$.

A more general approach for computation of the Thurston-Bennequin
number associated to a normal unibranch Q-singularity was partly
announced in \cite{Fin2} using a resolution of the singularity.
Below, Theorem \ref{calcul}, Section \ref{sec:6}, offers a method
to compute the Thurston-Bennequin number of a real normal
unibranch Q-singularity using a very good resolution graph of the
singularity. In the case the real part of the ambient real
algebraic surface is oriented, Theorem \ref{calcul} gives Equation
\ref{selfint}.

\subsection{Resolution graphs} \label{sec:2}
Let $\Cc X$ be a normal complex surface and
$\rho:\Cc\widetilde{X}\rightarrow\Cc X$ be a resolution of $\Cc
X$. The preimage $E=\rho^{-1}(\sing(\Cc
X))\subset\Cc\widetilde{X}$ is called the exceptional set and each
irreducible curve in the exceptional set is called an exceptional
curve of the resolution $\rho$. Let the number of exceptional
curves in $E$ be $s$; then the exceptional curves of $E$ will be
denoted by  $E_1,\ldots,E_s$.

The resolution $\rho$ is said to be {\it  good} if the exceptional
curves are smoothly embedded and any two intersecting exceptional
curves intersect at normal crossings. A {\it very good} resolution
is a good resolution with any pair of exceptional curves having at
most one common point. A good and a very good resolution always
exist for any analytic surface (cf. \cite{Lau} or \cite{Dur1}).

The dual graph of a resolution of a surface singularity is
constructed by representing each exceptional curve by a vertex and
connecting two vertices by an edge if the corresponding
exceptional curves intersect. It is well known that the dual graph
$\Gamma$ of a good resolution of a normal surface singularity is
connected. The $s$ vertices $e_1,\ldots,e_s$ of $\Gamma$
correspond to $E_1,\ldots,E_s$ respectively. Each vertex $e_i$
($1\leq i\leq s$) is endowed with a pair of integers: the self
intersection $n_i$ and the genus $g_i$ of $E_i$. Two vertices
$e_i$ and $e_j$ are connected by $n_{ij}$ edges where $n_{ij}$ is
the number of points in $E_i\cap E_j$.

From now on, let us assume that $\Gamma$ is the dual graph
of a very good resolution $\rho:\Cc\widetilde{X}\rightarrow\Cc X$
and also that $\Gamma$ is a tree.

\begin{figure}[h]
\begin{center}

{\tiny a rupture vertex}

\epsfig{file=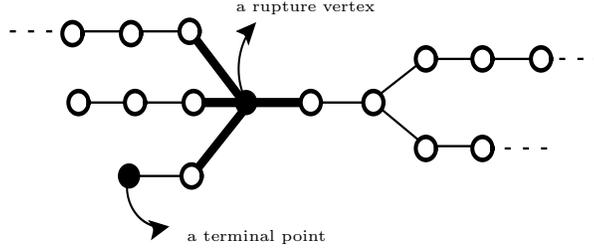,height=3cm,width=8cm}

\vspace{-0.45cm}

{\tiny a terminal point} \hspace{1.2cm}

\caption{A terminal point, a rupture vertex and its
star}\label{rupstar}
\end{center}
\end{figure}

If a vertex $e$ of $\Gamma$ is connected to more than two
vertices, then $e$ is called a {\it rupture vertex}. A vertex
adjacent to one and only one vertex of $\Gamma$ is called a {\it
terminal vertex}. The {\it star} of a subgraph $\sigma$ of
$\Gamma$ is a subgraph of $\Gamma$ consisting of $\sigma$ and all
the edges connected to any vertex of $\sigma$. The star of
$\sigma$ is denoted by $\sitare(\sigma)$ (see Figure
\ref{rupstar}).

Let $e$ be a vertex of $\Gamma$. Let us denote  by $\Gamma_e$
the subgraph $\Gamma-\sitare(e)$
of $\Gamma$. Since $\Gamma$ is a tree, the graph $\Gamma_e$ is
disconnected if $e$ is not a terminal vertex. Each connected component of
$\Gamma_e$ is called an {\it arm} of $e$.
The number of arms of $e$ is equal to the number of edges connected to $e$.

Let $e_1,\ldots,e_k$ be a set of vertices of a graph $\Gamma$. The
{\it subgraph spanned by the vertices $e_1,\ldots,e_k$} is the
subgraph of $\Gamma$ consisting of the vertices $e_1,\ldots,e_k$
and all the edges of $\Gamma$ connecting any pair of vertices
$e_i$ and $e_j$ with $1\leq i,j\leq k$.

The complex conjugation $\conj$ of $\Cc\widetilde{X}$ induced by
that of $\Cc X$ defines a map on the vertices of $\Gamma$ as the
following: \[\conj(e_i)=\bar{e}_i=\left\{\begin{array}{ll}
e_i,& \mbox{if $E_i$ is real}\\
e_j,& \mbox{if } \conj(E_i)=E_j.
\end{array}\right. \]
This map is 1-1 and onto since each exceptional curve $E_i$ is
either real or it does not contain any real point at all. In fact,
any normal singularity of a complex surface locally can be thought
of as a branched cover of $\Cc^2$ branched along a singular curve
\cite{Lau}. Hence a resolution of the singularity can be
constructed (as in Section \ref{sec:4}) by taking the branched
cover of a resolution of $\Cc^2$ and then by normalizing. Each
exceptional curve in the resolution of $\Cc^2$ can be made real so
that, in the cover, the preimage of each exceptional curve is
real. This said, the claim follows (for this fact see \cite{Fin3}
or for real resolutions see \cite{Sil}).

Now, if an exceptional curve is real,  the corresponding vertex
will be called real; otherwise it will be called an imaginary
vertex. Note that the map conj on the vertices of a resolution
graph is a symmetry of the graph.

Let us denote by $\Gamma_{\Rr}$ and $\Gamma_{\Cc}$ the subgraphs spanned by
the real and respectively imaginary vertices of $\Gamma$.
Furthermore, let ${\mathcal C}_e$ denote
the set of imaginary arms in $\Gamma_e$, that is, the arms of $e$
on which every vertex is imaginary.

Let $e$ be a rupture vertex of $\Gamma$. We will assign a weight
$n^{\sigma}$ to each imaginary arm $\sigma$ in ${\mathcal C}_e$
and a number $n'_e$ to vertex $e$. First assume that the arm
$\sigma\in{\mathcal C}_e$ is a {\it bamboo}, i.e. contains no
rupture vertex. Let there be $k$ vertices $e_1,\ldots,e_k$ on
$\sigma$ such that $e_i$ is closer than $e_j$ to $e$ if $i<j$
$(1\leq i,j\leq k)$. Then the weight of arm $\sigma$ is defined to
be the continued fraction:
\[n^{\sigma}={n_1-\frac{1}{n_2-\frac{1}{\cdots-\frac{1}{n_k}}}}.\]
If the imaginary arms of vertex $e$ are all bamboos, then define:
\[n'_e=n_e-\sum_{\sigma\in{\mathcal C}_e} \frac{1}{n^{\sigma}}.\]
Now, let an imaginary arm $\sigma\in{\mathcal C}_e$ consist of
vertices $e_1,\ldots,e_k$, some of which may be rupture vertices.
For each rupture vertex $e_i$ on $\sigma$ compute $n'_{e_i}$
recursively; furthermore, put $n'_i=n_i$ if $e_i$ is a vertex on
$\sigma$ that is not a rupture vertex. For a short hand notation,
we put $n'_i=n'_{e_i}$. Then define
\[n^{\sigma}={n'_1-\frac{1}{n'_2-\frac{1}{\cdots-\frac{1}{n'_k}}}}.\]

\subsection{Embedded resolution graphs} \label{sec:3}
Let us consider the reduced real algebraic curve $\Cc
A=\{f(x,y)=0\}$ in $\Cc^2$ with a singularity at 0 and with no
other singular point. One can resolve $\Cc A$ by a series of
blow-ups of $\Cc^2$. Topologically the blown-up surface is
$\Cc^2\#t\overline{\Cc P^2}$ where $t$ is the number of
blowing-ups. Let us denote the resolution by
$\rho':\Cc\widetilde{Y}\rightarrow \Cc^2$ where
$\Cc\widetilde{Y}=\Cc^2\#t\overline{\Cc P^2}$. The preimage
$(\rho')^{-1}(\Cc A)$ as a divisor is called the {\it total
transform} of $\Cc A$. The closure of the preimage
$(\rho')^{-1}(\Cc A-\{0\})$ as a divisor is called the {\it proper
transform} of $\Cc A$, denoted by $\Cc\widetilde{A}$. Each
component of $\Cc\widetilde{A}$ is smooth and can be made to
transversally intersect the exceptional curves by sufficiently
many blow-ups at the points of non-transversal intersection.

Let us furnish  the resolution graph corresponding to the
resolution $\rho'$ with the following further data. Let the number
of exceptional curves of $\rho'$ be $t$. Each component of the
proper transform $\Cc\widetilde{A}$ which intersects an
exceptional curve $E_i$ ($1\leq i\leq t$) is specified in the
resolution graph by an arrow attached to the vertex $e_i$. For
instance, if $f(x,y)$ is irreducible, then there exists one arrow
on the graph. Furthermore, each vertex $e_i$ is assigned a
multiplicity $m_i\in\Zz$; each arrow is assigned a multiplicity 1.
Such an `extended' graph is called an {\it embedded resolution
graph} and is denoted by $\Gamma_f$. It is easy to see that
$\Gamma_f$ is always a tree with arrows.

We extend the definition for a rupture vertex on an embedded
resolution graph as follows: A vertex of $\Gamma_f$ which is
adjacent to more than two vertices or two vertices and an arrow is
called a {\it rupture vertex} of $\Gamma_f$.

The multiplicity $m_i$ of an exceptional curve $E_i$ is the
vanishing order of the map
$f\circ\rho':\widetilde{Y}\rightarrow\Cc$ on $E_i$. In other
words, in the chart of $\widetilde{Y}$ containing the intersection
point of $E_i$ and $E_j$, the neighbourhood of the intersection
can be viewed as $\{(x,y)\in\Cc^2|x^{m_i}y^{m_j}=0\}$ with an
appropriate choice of local coordinates $(x,y)$.

It is well known that the following identity is satisfied on
$\Gamma_f$ for every $k$, $1\leq k\leq t$ (cf. e.g. \cite{GoS}, p.
251):
\begin{equation}n_k m_k+\sum_{\stackrel{1\leq i\leq t}{i\neq k, n_{ik}\neq 0}} m_i =0.
\label{mini} \end{equation}

\begin{figure}[h]
\begin{center}

\hspace{0.3cm} $m_i:\;\; n$ \hspace{1.1cm} $q_l n$ \hspace{1.7cm}
$\cdots$ \hspace{1.4cm} $q_{l-1}m$ \hspace{0.8cm} $q_l
n\!+\!r_l=m$

\epsfig{file=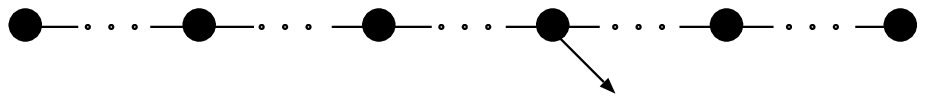,height=1cm,width=8.5cm}

\vspace{-0.7cm}

$\underbrace{\mbox{ \hspace{1.8cm}}}_{\displaystyle q_l \mbox{
times}}$ \hspace{1.2cm} $\underbrace{\mbox{
\hspace{1.8cm}}}_{\displaystyle q_0 \mbox{  times}}$
\hspace{1.1cm} $\underbrace{\mbox{ \hspace{1.8cm}}}_{\displaystyle
q_{l-1} \mbox{  times}}$

\caption{} \label{E00}
\end{center}
\end{figure}

As a special case, it is going to be useful to understand the
algorithmic construction of the embedded resolution graph for
$f(x,y)=x^m+y^n$, $\gcd(m,n)=1$. Assume that $m>n$. One can
immediately observe that the blow-up sequence associated to the
singularity of $f(x,y)=0$ at 0 is directly related with the
Euclidean algorithm for $m$ and $n$:
\[m=q_l n+r_l;\]
\[n=q_{l-1} r_l+r_{l-1};\]
\[r_l=q_{l-2} r_{l-1}+r_{l-2};\]
\[\cdots\]
\[r_3=q_1 r_2+1;\]
\[r_2=q_0\cdot 1\] where $r_l<n$, $r_i\in\Zz^+$ and $r_{i-1}<r_i$ for all $i$ with $1<i<l$.
The line containing $q_i$ above corresponds to $q_i$ consecutive
blow-ups so that there are $\sum_{i=0}^l q_i$ vertices of
$\Gamma_f$. Let us index the vertices of $\Gamma_f$ in the
following way: $e_{q_i,j}$ with $0\leq i\leq l$ and $1\leq j\leq
q_i$ denotes the vertex created by $j^{\rm th}$ blow-up within the
$q_i$ consecutive blow-ups. The unique arrow is attached to the
vertex $e_{q_0,q_0}$. The graph $\Gamma_f$, the vertices
$e_{q_i,j}$ and their multiplicities are shown in Figure
\ref{E00}.

The unique rupture vertex $e_{q_0,q_0}$ of $\Gamma_f$ will be
denoted by $e_0$. The arm of $e_0$ with the terminal vertex of
multiplicity $n$ (respectively $m$) will be called an {\it
($n$)-arm} (respectively {\it ($m$)-arm}) of $\Gamma_f$.

\section{Resolution graphs of Brieskorn double points}
\label{sec:4}
Let $f(x,y)=x^m+y^n$, $\Cc X$ be the double branched cover of
$\Cc^2$ determined by $f(x,y)\pm z^2=0$, the branching locus being
$\Cc A=\{f(x,y)=0\}$. Note that the Brieskorn double surface $\Cc
X$ has a normal Q-singularity at 0 (see e.g. \cite{Sav}, Chapter
1).

We are going to consider the commutative diagram:
\[
\begin{array}{rcl}
\Cc\widetilde{X}&\stackrel{\textstyle\rho}{\longrightarrow} &\Cc
X\\ \widetilde{p}\Big\downarrow &&\Big\downarrow p\\
\Cc\widehat{A}\subset\Cc\widetilde{Y}&\stackrel{\textstyle\rho'}{\longrightarrow}
&\Cc^2\supset\Cc A
\end{array}
\]
where $p:(x,y,z)\mapsto (x,y)$ is the projection,
$\rho:\Cc\widetilde{X}\rightarrow\Cc X$ is a very good resolution
and $\widetilde{p},\rho'$ are defined as follows.

The non-singular surface $\Cc\widetilde{X}$ can be constructed in
two steps. This construction, which is sketched below, can be
traced, for example, in \cite{HKK}, Section 0; \cite{GoS} Section
7.2 or \cite{Nem}. At the first step, the embedded resolution
graph $\Gamma_f$ is constructed. Furthermore, if two exceptional
curves with odd multiplicities intersect, they are made disjoint
by an additional blow-up at the point of intersection. Denote the
graph obtained in that way by $\Gamma'_f$ and the corresponding
resolution by $\rho':\Cc\widetilde{Y}=\Cc^2 \#s\overline{\Cc
P}^2\rightarrow\Cc^2$ where $s$ is the number of blow-ups employed
so far.

Next, let $\Cc\widehat{A}\subset\Cc\widetilde{Y}$ be the union of
the proper transform $\Cc\widetilde{A}=(\rho')^{-1}(\Cc A)$
and the exceptional curves in $\Cc\widetilde{Y}$ with odd
multiplicity. The double covering of $\Cc\widetilde{Y}$ branched along
$\Cc\widehat{A}$ with an appropriate choice of real structure gives
$\Cc\widetilde{X}$. The covering map is denoted by $\widetilde{p}$.

If $E_i$ is an exceptional curve in $\Cc\widetilde{Y}$ with
odd multiplicity, then
the preimage $\widetilde{p}^{-1}(E_i)$ is connected
and has self intersection $\frac{n_i}{2}$.
An exceptional curve $E_i$ of even multiplicity intersects either with two
vertices of odd multiplicity or with no vertices of odd multiplicity.
In case of former, $\widetilde{p}^{-1}(E_i)$ is connected and has
self intersection $2n_i$. In the latter case, $\widetilde{p}^{-1}(E_i)$ is
constituted of two exceptional curves, each with self-intersection $n_i$.

The graph constructed in this way is a good resolution graph and can be
changed into a minimal good resolution graph by blowing-down the
$(-1)$-curves that do not intersect more than two exceptional curves.
We are going to denote by $\Gamma(m,n)$ the resolution graph
corresponding to $\rho$.

Let the graph $\Gamma(m,n)$ have $r$ vertices. and let $E^i$
(respectively $e^i$) ($1\leq i\leq s$) denote the exceptional
curves of $\rho$ (respectively the corresponding vertices of
$\Gamma(m,n)$). The covering map
$\widetilde{p}:\Cc\widetilde{X}\rightarrow\Cc\widetilde{Y}$
induces a map $\widetilde{p}:\Gamma(m,n)\rightarrow\Gamma_f$,
sending $e^i$ to $e_j$ if $E^i$ lies in
$\widetilde{\rho}^{-1}(E_j)$.

The following facts on $\Gamma_f$ and $\Gamma(m,n)$ are direct
consequences of the algorithm summarised above and at the end of
previous section.

\begin{lemma}
Let $f(x,y)=x^m+y^n$ $(\gcd(m,n)=1)$, $e_0$ be the rupture vertex
of $\Gamma_f$ and $\widetilde{p}:\Gamma(m,n)\rightarrow\Gamma_f$
be the covering map induced by $\widetilde{p}:\Cc X=\{x^m+y^n\pm
z^2=0\}\rightarrow\Cc^2 \#t\overline{\Cc P}^2$. Then:

(i) Multiplicity of $e_0=m_{e_0}=mn$.

(ii) $\widetilde{p}^{-1}(e_0)$ is one vertex. Let us denote it by $e^0\in\Gamma(m,n)$.

(iii) Let $\sigma$ be an arm of $e_0$. Each connected subgraph of
$\Gamma(m,n)$ in the preimage $\widetilde{p}^{-1}(\sigma)$ is an arm of $e^0$.

(iv) There are $\gcd(m,2)$ arms of $e^0$ in the preimage of an
($n$)-arm of $\Gamma_f$. Analogously, there are $\gcd(n,2)$ arms
of $e^0$ in the preimage of an ($m$)-arm.

(v) Number of arms of $e^0$ is 3 and each arm of $e^0$ is a bamboo.
\end{lemma}

The fact (iii) above allows the following definition:
each arm of $\Gamma(m,n)$ in the preimage
of an ($n$)-arm of $\Gamma_f$ will be called an {\it ($n$)-arm} of $\Gamma(m,n)$.
Analogous definition is made for an {\it ($m$)-arm} of $\Gamma(m,n)$.

Note that any real algebraic surface defined as a double cover of
$\Cc^2$ branched along the reduced curve $\{f(x,y)=0\}$ can be
assigned two complex conjugations (see e.g. \cite{DeK}, Section
2.2.4). For the case of Brieskorn double surfaces when
$f(x,y)=x^m+y^n$, let us denote the complex conjugations by
$\conj_\pm$ and the corresponding complex surfaces $\{x^m+y^n\pm
z^2=0\}$ by $\Cc X_{\pm}$. First, we give a proof of the following
fact pointed out in \cite{Fin2}:

\begin{proposition}
\label{prop0} Let $\Gamma_{\Rr}^{\pm}(m,n)$ be the subgraph
spanned by the vertices of $\Gamma(m,n)$ fixed under the action of
$\conj_\pm$. If $m$ and $n$ are odd, then
\[\Gamma_{\Rr}^{\pm}(m,n)=\Gamma(m,n).\] If $m$ is odd and $n$ is even, then \[
\Gamma_{\Rr}^-(m,n)=\Gamma(m,n)\] and $\Gamma_{\Rr}^+(m,n)$ is the
subgraph of $\Gamma(m,n)$ spanned by the rupture
vertex and the the ($n$)-arm.
\end{proposition}

\begin{proof}
If $m$ and $n$ are both odd, then the resolution graph
$\Gamma(m,n)$ has no symmetry but the identity. Since complex
conjugations correspond to symmetries of the resolution graph,
both choices of complex conjugation coincide and all exceptional
curves are fixed under conjugation of either choice so that the
first assertion follows.

Now let us assume  that $m$ is odd and $n$ is even. Let
$f(x,y)=x^m+y^n$ and let $e_0$ be the rupture vertex of
$\Gamma_f$. The multiplicity $m_0$ of $E_0$ is $mn$ and the self
intersection is $-1$. The embedded resolution graph $\Gamma_f$
around the vertex $e_0$ is shown in Figure \ref{E0}(a). The two
vertices $e_1$ and $e_2$ adjacent to $e_0$ have multiplicities of
different parities by Equation \ref{mini}, Section \ref{sec:3}.
Let us assume that $m_1$ is odd and $m_2$ is even. In a properly
chosen chart of $\Cc^2\#t\overline{\Cc P}^2$ around the point of
intersection $E_0\cap E_1$, the configuration of the exceptional
curves $E_0$, $E_1$, $E_2$ and the proper transform
$\Cc\widetilde{A}$ is shown in Figure \ref{E0}(b). The curve
$E_1\cup\Cc\widetilde{A}$ is in the branching locus and is
expressed by $x^{mn}(x+1)=0$ in this chart.

\begin{figure}[h]
\begin{center}

\hspace{0.6cm} $e_1$ \hspace{0.7cm} $e_0$ \hspace{0.7cm} $e_2$
\hspace{3.0cm} $\Cc\widetilde{A}$ \hspace{0.3cm} $E_1$
\hspace{0.7cm} $E_2$

\vspace{0.2cm}

\hspace{2.8cm} $E_0$

\vspace{0.1cm}

\hspace{3.4cm} $-1$

\vspace{0.1cm}

\hspace{4.1cm} $0$

\vspace{-1.6cm}

\epsfig{file=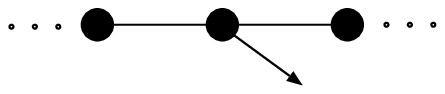,height=1cm,width=4.3cm} \hspace{2cm}
\epsfig{file=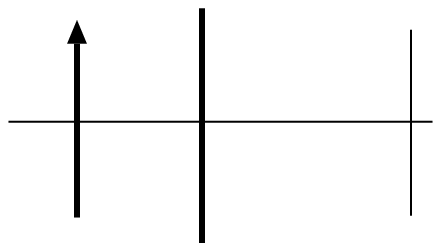,height=1.6cm,width=3cm}

\caption{}
\label{E0}
\end{center}
\end{figure}

Suppose that the double covering $\Cc X$ of $\Cc^2$ branched along $\Cc A$ is determined
by $\{x^m+y^n-z^2=0\}$. Then above the subset
${\mathcal S}=\{(x,y)|-1<x<0\}$ of the chart in Figure \ref{E0}(b)
lie the imaginary points of $\Cc X$ and above the points outside ${\mathcal S}$
on the chart lie the real points of $\Cc X$. In particular, each
connected component of the preimage $\widetilde{p}^{-1}(E_2)$
is real so that all the exceptional curves on its arm have real connected
components on their preimages.
Also $\widetilde{p}^{-1}(\Cc\widetilde{A})$, $\widetilde{p}^{-1}(E_0)$
and $\widetilde{p}^{-1}(E_1)$ are real in $\Cc X$.

Alternatively, suppose that $\Cc X=\{x^m+y^n+z^2=0\}$. It follows with an
analogous argument to the previous paragraph that
this time each connected component of the preimage $\widetilde{p}^{-1}(E_2)$ contains
no real point because of the sign change.
\end{proof}

\begin{proposition}
\label{prop1} Let \[n'=\left\{\begin{array}{ll}
n+2m & m\;\; \odd\\n+m & m\;\; \even.\end{array}\right.\]
For $f(x,y)=x^m+y^n$ and
$f'(x,y)=x^m+y^{n'}$, the resolution graph $\Gamma_{f'}$ differs
from $\Gamma_f$ by a pair of terminal vertices. The outermost
vertex has multiplicity $m$ and self intersection $-2$ while the
adjacent vertex has multiplicity $2m$ (see Figure
\ref{2nodes}). Furthermore, the resolution graph $\Gamma(m,n')$
differs from $\Gamma(m,n)$ by a terminal vertex on each ($n$)-arm of
$\Gamma(m,n')$.
\end{proposition}

\begin{proof}
Let us denote the surfaces
$\{f(x,y)+z^2=0\}$ and $\{f'(x,y)+z^2=0\}$ by $X$ and $X'$ and the
branching curves $\{f(x,y)=0\}$ and $\{f'(x,y)=0\}$ by $A$ and
$A'$ respectively.

\begin{figure}[h]
\begin{center}

$e_1\;\;\;\;\;\;\;\;e_2$ \hspace{1.2cm}

\epsfig{file=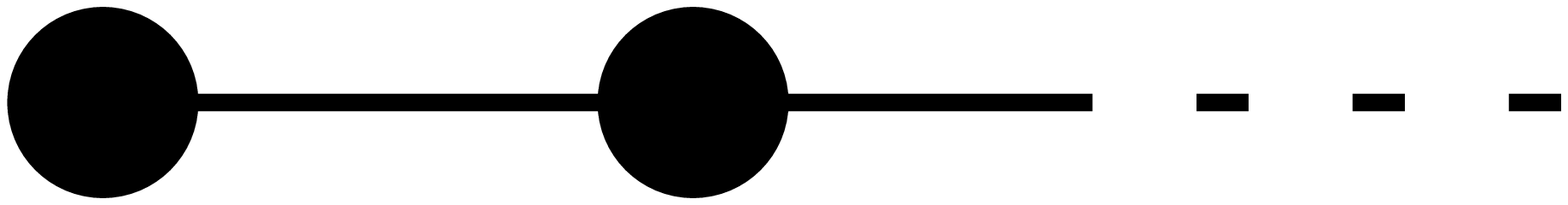,height=0.35cm,width=2.8cm}

${\displaystyle {m\atop -2}\;\;\;\;\;{2m\atop n_2}}$
\hspace{1.2cm}

\caption{$\Gamma_{f'}$ has two extra vertices for odd $m$}\label{2nodes}
\end{center}
\end{figure}

First suppose that $m$ is odd.
Let us blow-up $\Cc^2$ once at 0 to obtain a proper transform
$A'_1\subset\Cc^2\#\overline{\Cc P}^2$ of $A'$ and an exceptional
curve $E_1$ with self intersection $n_1=-1$ and multiplicity $m_1=m$.
Then let us blow-up the
$\Cc^2$-chart of $\Cc^2\#\overline{\Cc P}^2$ containing the
singular point of $A'_1$. The exceptional set contains the exceptional curve
$E_1$ with $n_1=-2$ transversally intersecting an exceptional
curve $E_2$ with $n_1=-1$, $m_2=2m$.
Furthermore the proper transform $A'_2$ of $A'_1$ has a singularity
at $A'_2\cap E_2$ given locally by $f(x,y)=x^m+y^n=0$, identical to the
curve $A$.

Further blow-ups to resolve the singularity of $A'_2$
increase the self intersection of $E_2$ at least by 1 and produce
a proper transform $A'_{t}$ of $A'$
and the exceptional curves $E_3,\ldots,E_{t}$. The curves $A',
E_3,\ldots,E_t$ are transverse to $E_2$
and constitute the graph $\Gamma_f$ which is a subgraph of $\Gamma_{f'}$.
Therefore the embedded resolution graph $\Gamma_{f'}$ associated to $f'$
differs from $\Gamma_f$ by the two vertices $e_1$
and $e_2$ appended to one end of $\Gamma_f$.

Now consider the very good resolution
$\rho':\widetilde{X'}\rightarrow\Cc^2\#t\overline{\Cc P^2}$
associated to $\Gamma_{f'}$. The
preimages $(\rho')^{-1}(E_1)$ and $(\rho')^{-1}(E_2)$ are rational
curves  with self intersections $-1$ and $-2n_2\leq -4$
respectively. Blowing $(\rho')^{-1}(E_1)$ down, we end up with a
vertex on  $\Gamma(m,n')$
with self intersection at most $-3$ (before further possible
blow-downs of $-1$ curves). This vertex is appended to the
($n$)-arm of $\Gamma(m,n)$.

For the case of even $m$, the proof is analogous except that in
the beginning one needs to blow-up only once to obtain a new
vertex in $\Gamma(m,n')$ representing an exceptional curve of
multiplicity $m$ and self intersection $-k$. Also note that,
different from the odd case, there are two ($n$)-arms on $\Gamma_f$
and two ($n'$)-arms on $\Gamma_{f'}$.
\end{proof}

\begin{corollary}
\label{cor1}
$\Gamma(m,n+4m/(m,2))$ differs from $\Gamma(m,n)$ by two extra
vertices at the end of each of its ($n$)-arms.
\end{corollary}

\begin{proposition}
\label{prop2} The self intersection of the terminal vertex of each
($n$)-arm of the graph $\Gamma(m,n+4m/(m,2))$ is $-2$.
\end{proposition}

\begin{proof}
Assume first that $m$ is odd. Let $n'=n+4m$ and
$f'(x,y)=x^m+y^{n'}$. One has the terminal
part of an $(n')$-arm of $\Gamma_{f'}$ as shown in
Figure \ref{4nodes}(a). To see this one needs to
apply four consequent blow-ups as described in the proof of
Proposition \ref{prop1}. Then it is straightforward to verify
that the terminal part of an ($n'$)-arm of $\Gamma(m,n')$ is as in
Figure \ref{4nodes}(b). Blowing down the $(-1)$-curves, the proof
follows.

\begin{figure}[h]
\begin{center}
$\!e_1'\;\;\;\;\;\;\;\;\,e_1\;\;\;\;\;\;\;\;\,e_2'\;\;\;\;\;\;\;\;e_2$\hspace{2.2cm}
$e_1'\;\;\;\;\;\;\;\,e_1\;\;\;\;\;\;\;\,e_2'\;\;\;\;\;\;\;\,e_2$\hspace{.5cm}

\hspace{0.5cm} \epsfig{file=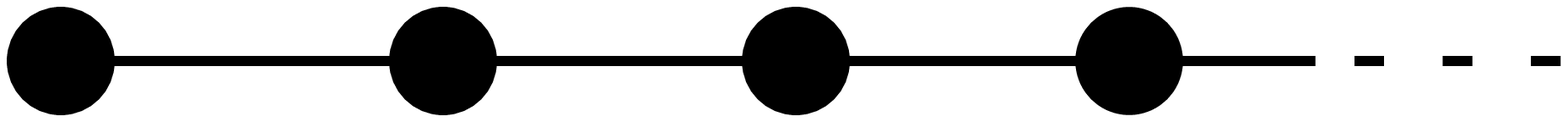,height=0.4cm,width=4.5cm}
\hspace{1.cm} \epsfig{file=fig5.eps,height=0.4cm,width=4.5cm}

\hspace{.1cm} ${\displaystyle {m\atop -2}\;\;\;\;{2m\atop
-2}\;\;\;\;\;\; {3m\atop -2}\;\;\;\;\;\;{4m\atop n_2}}$
\hspace{1.9cm} $-1\;\;\;\;\;-4\;\;\;\;\;-1\;\;\;\;\;\;\;\;2n_2$
\hspace{.5cm}

\vspace{.5cm}

(a) \hspace{7cm} (b)

\caption{Resolution graphs (a) $\Gamma_{f'}$, (b)
$\Gamma(m,n+4m)$ for odd m}\label{4nodes}
\end{center}
\end{figure}

The case for even $m$ is similar.
\end{proof}

\section{The characteristic set of a resolution graph and the
first Chern class} \label{sec:5}
Let us consider a very good resolution
$\rho:\Cc\widetilde{U}\rightarrow\Cc U$ of a compact cone-like
neighbourhood $\Cc U$ of a normal singularity and denote
$\Cc\widetilde{M}=\partial\Cc\widetilde{U}$. As in \cite{Dur2},
let $c$ be the unique class in
$H^2(\Cc\widetilde{U},\Cc\widetilde{M};\Qq) \equiv
H^2(\Cc\widetilde{U},\Cc\widetilde{M})\otimes\Qq$ such that
$j^*(c)=c_1(\Cc\widetilde{U}) \in H^2(\Cc\widetilde{U};\Zz)$ where
$j:\Cc\widetilde{U}\rightarrow
(\Cc\widetilde{U},\Cc\widetilde{M})$ is the inclusion map. The
canonical class $K\in H_2(\Cc\widetilde{U};\Qq)$ is defined to be
the Lefschetz dual $LD(-c)$ and can be expressed as $\sum_{1\leq
i\leq s} a_i E_i$ where $a_i\in \Qq$ and $s$ is the number of
exceptional curves of the resolution.

If the singularity above is numerically Gorenstein, i.e.
$K$ is integral, then $K$ is integral dual to the second Stiefel-Whitney class
$w=w_2(\Cc\widetilde{U},\Cc\widetilde{M};\Zz)$. Therefore, we can write
$w=\sum_{1\leq i\leq s} \hat{a}_i E_i$ where $\hat{a}_i=(a_i\!\mod\!2)$ and
$E_i$ is identified with its dual in
$H^2(\Cc\widetilde{U},\Cc\widetilde{M};\Zz)$. Note that $w$ satisfies the Wu formula:
\begin{equation} w\cdot E_i\equiv E_i\cdot E_i \mbox{ (mod 2) for all } i, 1\leq i\leq s
\label{defw} \end{equation}
and furthermore for $1\leq i,j\leq s$,
the matrix of the intersection form $Q_{\Cc\widetilde{U}}$ on $\Cc\widetilde{U}$ has the entries
\[E_i\cdot E_j=\left\{
\begin{array}{ll}n_{ij},&i\neq j\\n_i,&i=j\end
{array}\right.\] in the basis $E_1,\ldots,E_s$.

For a normal, numerically Gorenstein
Q-singularity, let $\Gamma$ be a very good resolution graph with $s$
vertices and $w$ be expressed as a sum as above.
The {\it characteristic set} ${\mathcal W}$ of $\Gamma$ is defined as the following set
of vertices of $\Gamma$:
\[{\mathcal W}=\{e_i\in\Gamma|\hat{a}_i=1\}.\]

It is algorithmic to determine the characteristic set of a very good
resolution graph. Nevertheless, the relation between $c_1(\Cc\widetilde{U})$ and
${\mathcal W}$ may lead to some general properties.

Now let $p:\Cc X\rightarrow\Cc ^2$ be a double covering branched along
an irreducible singular curve $\Cc A=\{f(x,y)=0\}$ that has an isolated singularity at 0.
Note that the singularity of $\Cc X$ is normal and numerically Gorenstein (see e.g. \cite{Dur2}).

As in Section \ref{sec:4}, consider the commutative diagram:
\[
\begin{array}{rcl}
\Cc\widetilde{U}&\stackrel{\textstyle\rho}
{\longrightarrow} &\Cc U\\
\widetilde{p}\Big\downarrow &&\Big\downarrow p\\
\Cc\widetilde{V}&\stackrel{\textstyle\rho'}{\longrightarrow}
&B^4\subset\Cc^2
\end{array}
\]
where $B^4$ is a sufficiently small ball around 0 so that $\Cc U$ is cone-like.
Let $\widetilde{A}$ denote the proper transform of $\Cc A\cap B^4$ under $\rho'$
and $\widehat{A}\subset\Cc\widetilde{V}$ the branching locus of the
double covering  $\widetilde{p}$; $E^i\subset\Cc\widetilde{U}$
($1\leq i\leq s$) denote the exceptional curves
of the very good resolution $\rho$; $E_i\subset\Cc\widetilde{V}$
($1\leq i\leq t$) denote the exceptional curves of the resolution $\rho'$.
Denote the multiplicity of $E_i$ by $m_i$.
We can express \[c_1(\Cc\widetilde{U})=\sum_{i=1}^s a_i E^i,\;\;\;\;\;
c_1(\Cc\widetilde{V})=\sum_{i=1}^s b_i E_i.\]
Above and what follows below, $E^i\in H_2(\Cc\widetilde{U})$
is always identified with the image $j_U^*(LD(E_i))$ of its Lefschetz dual $LD(E_i)\in
H^2(\Cc\widetilde{U},\partial\Cc\widetilde{U})$
under the homomorphism $j_U^*$ induced by the inclusion
$j_U:\Cc\widetilde{U}\rightarrow(\Cc\widetilde{U},\partial\Cc\widetilde{U})$
and $E_i$ is identified with the image $j_V^*(LD(E_i))$
under the inclusion
$j_V:\Cc\widetilde{V}\rightarrow(\Cc\widetilde{V},\partial\Cc\widetilde{V})$.

\begin{proposition}
\label{lem3} Assume that $E_j=\widetilde{p}(E^i)$ for some $i$ and $j$ with
$1\leq i\leq s$ and $1\leq j\leq t$. If the
multiplicity $m_j$ of $E_j$ is odd then the coefficient $a_i$ in
$c_1(\Cc\widetilde{U})$ is even.
\end{proposition}

\begin{proof}
First note that $c_1(\Cc\widetilde{U})=(\widetilde{p})^*
(c_1(\Cc\widetilde{V})-\frac{1}{2}[\hat{A}])$
(see e.g. \cite{Bea}).  The class
$c_1(\Cc\widetilde{V})-\frac{1}{2}[\hat{A}]$ can be expressed
as a sum $\sum_{i=1}^t d_i
E_i$ in $H^2(\Cc\widetilde{V})$ with $d_i\in\Zz$. Then,
\[c_1(\Cc\widetilde{U})=(\widetilde{p})^*(\sum_{i=1}^t d_i E_i)
={\displaystyle \sum_{1\leq i\leq t\atop m_i \odd}} 2d_i
(\widetilde{p})^*(E_i)+{\displaystyle \sum_{1\leq i\leq t\atop m_i \even}} d_i
(\widetilde{p})^*(E_i)\] which shows that all of the exceptional
curves with odd multiplicities have even coefficients in
$c_1(\Cc\widetilde{U})$.
\end{proof}

\begin{proposition}
\label{mel}
Assume that $E_j=\widetilde{p}(E^i)$ for some $i$ and $j$ with
$1\leq i\leq s$ and $1\leq j\leq t$. Also let $(\widetilde{p})^*(\frac{1}{2}[\widehat{A}])
=\sum_{k=1}^s d_k E^k$, $d_k\in\Zz$. Then $d_i$ is odd if and only if
the multiplicity $m_j\equiv 2\mod\,4$.
\end{proposition}

\begin{proof}
First let us note:
\begin{lemma} \hspace{1cm} ${\displaystyle [\widetilde{A}]=-\sum_{l=1}^t m_l E_l.}$
\end{lemma}
\begin{proof}
The proper
transform $\widetilde{A}$ intersects just one exceptional curve,
say $E_r$. To prove the claim, it is enough to check that
\[Q_{\Cc\widetilde{V}}\bigg((-\sum_{l=1}^t m_l E_l),E_r\bigg)=\left\{
\begin{array}{ll}1,&r=l\\0,&r\neq l.\end{array}\right.\] It is then
straightforward to see that this identity is in fact a
reformulation of Equation \ref{mini}, Section \ref{sec:3}.
\end{proof}

Since \[[\widehat{A}]=[\widetilde{A}]+\sum_{m_i \odd}E_i,\]
it follows from the lemma that
\[
\begin{array}{ll} [\widehat{A}]&=-{\displaystyle \sum_{i=1}^t}
m_i E_i+{\displaystyle \sum_{1\leq i\leq t\atop m_i \odd}}E_i\\ &=-{\displaystyle
\sum_{1\leq i\leq t\atop m_i \even}} m_i E_i+{\displaystyle \sum_{1\leq i\leq t\atop
m_i \odd}} (-m_i+1) E_i.
\end{array}\]

Considering the behaviour of the transfer map from
$H^2(\Cc\widetilde{U};\Zz)$ to $H^2(\Cc\widetilde{V};\Zz)$  one can deduce that
\[(\widetilde{p})^*\Big(\frac{1}{2}[\widehat{A}]\Big)=\sum_{1\leq k\leq t\atop m_k \even}
\frac{m_k}{2} (\widetilde{p})^*\Big(E_k\Big)
+\sum_{1\leq k\leq t\atop m_k \odd}(-m_k+1)
(\widetilde{p})^*\Big(E_k\Big).\]
A coefficient on the right hand side corresponding to $E_j$ is odd if and only if
$m_j\equiv 2\mod\,4$.
\end{proof}

The last two propositions now prove
\begin{proposition}
Let $E_j=\widetilde{p}(E^i)$ for some $i$ and $j$ with
$1\leq i\leq s$ and $1\leq j\leq t$. Then $a_i$ is odd if and only if
one of the following two conditions holds: (i) $m_j\equiv 2\mod\,4$ and $b_j$ even;
(ii) $4|m_j$ and $b_j$ odd.
\label{tekcift}
\end{proposition}

\begin{proposition}
\label{lem4} Corresponding to the resolution $\rho':\Cc\widetilde{V}\rightarrow B^4$,
consider the blow-up sequence:
\[
\begin{array}{ccccccccccc}
\Cc\widetilde{V}=Y_t&\stackrel{\textstyle\rho_t}{\longrightarrow}&\ldots
&\stackrel{\textstyle\rho_{k+1}}{\longrightarrow}&Y_k&
\stackrel{\textstyle\rho_k}{\longrightarrow}&Y_{k-1}&
\stackrel{\textstyle\rho_{k-1}}{\longrightarrow}&\ldots&
\stackrel{\textstyle\rho_1}{\longrightarrow}&Y_0=B^4\subset\Cc^2\\
\cup&&&&\cup&&\cup&&&&\cup\\
\widetilde{A}&\longrightarrow&\ldots&\longrightarrow&A_k&\longrightarrow&A_{k-1}&
\longrightarrow&\ldots&\longrightarrow&A_0=\Cc A\cap B^4
\end{array}
\] where $Y_i$ ($1\leq i\leq t$) is obtained by blowing up $Y_{i-1}$
at the singular point of $A_{i-1}$, $\rho_i$ is the corresponding map
and $A_i$ is the proper transform of $A_{i-1}$. Let $E_k$ denote the
exceptional sphere $Cl(\rho_k^{-1}(A_{k-1})-A_k)$ in $Y_k$. Assume that
$c_1(Y_k)={\displaystyle \sum_{i=1}^k} b^{(k)}_i E_i$,
$b^{(k)}_i\in\Zz$. Then,
\[c_1(Y_{k+1})=c_1(Y_k)
+\Big(\!\!\!-1+\!\!\!\!\sum_{1\leq i\leq k \atop E_i\cap E_k\neq\emptyset}
\!\!b^{(k)}_i\Big) E_{k+1}.\]
\end{proposition}

\begin{proof} It is enough to show that $c_1(Y_{k+1})$ as expressed above
satisfies the adjunction formula $E_i\cdot
E_i+2=c_1(Y_{k+1})(E_i)$ for all $i$, $1\leq i\leq k+1$. Let
$n_i^{(k)}$ denote the self intersection of $E_i$ in $Y_k$. Then
the self intersection $n_i^{(k+1)}$ of $E_i$ in $Y_{k+1}$ is
decreased by 1 if $E_i\cap A_k\neq\emptyset$ in $Y_k$ and is
unchanged otherwise. Furthermore, $n_{k+1}^{(k+1)}=-1$ and by
Equation \ref{mini}, Section \ref{sec:3}, $m_{k+1}=1+\!\!\!\!\!\!
{\displaystyle \sum_{1\leq i\leq k \atop E_i\cap
E_{k+1}\neq\emptyset}}\!\!\!\!\! m_i$. It is now straightforward
to check that $c_1(Y_{k+1})$ as given in the form above satisfies
the adjunction formula.
\end{proof}

\begin{remark}
\label{lem5} Assume that $c_1(Y_k)=\sum_{i=1}^k d_i E_i$.
It follows from the previous proposition that \[c_1(Y_{k-1})=\sum_{i=1}^{k-1} d_i E_i.\]
\end{remark}

\begin{proposition} Let $\Cc A=\{f(x,y)=x^m+y^n=0\}\cap B^4$ with $\gcd(m,n)=1$.
Assume that $E_i$ ($1\leq i\leq t$) are as above and
that $c_1(\Cc\widetilde{V})=\sum_{i=1}^t b_i E_i$. Then
$b_t$ is odd if and only if both $m$ and $n$ are odd.
\label{b0cift}
\end{proposition}

\begin{proof} Let us adopt the notation of the remark at the end of
Section \ref{sec:3} so that $t=\sum_{k=0}^l q_k$. Using
Proposition \ref{lem4}, one can algorithmically calculate the
coefficients in the expression of $c_1(\Cc\widetilde{V})$. Let
$b_{q_i,j}$ be the coefficient corresponding to vertex
$e_{q_i,j}$. Then, $b_{q_l,1}=-1, b_{q_l,2}=-2\ldots
b_{q_l,q_l}=-q_l, b_{q_{l-1},1}=-q_l-1,
b_{q_{l-1},2}=-2q_l-2\ldots
b_{q_{l-1},q_{l-1}}=-q_{l-1}q_l-q_{l-1}\ldots$ Then it is
straightforward to verify by induction on $l$ that:
\[-b_t=-b_{q_0,q_0}=q_l n+q_{l-1} r_l+q_{l-2} r_{l-1}
+\ldots+q_1 r_2+q_0\cdot 1=m+n-1.\]
Hence the proof follows.
\end{proof}

\section{Computation of the invariant via resolution graphs} \label{sec:6}
I give the proof of the following fact which was stated in
\cite{Fin2}.

\begin{proposition}
Let $\Cc E$ be a real rational $(-m)$-curve $(m>0,m\in\Qq)$ in a
real algebraic Q-surface $\Cc X'$. Let $\Rr E$ be smooth and $\Rr
X'$ be nonorientable. Consider the projection $p:\Cc X'\rightarrow
X=\Cc X'/\Cc E$. Then the point $x=p(\Cc E)\in X$ is a
Q-singularity. Furthermore one can define a bilinear form $\tb_x$
on $H_1(M)$ (where $M=S_\epsilon\cap p(\Rr X')$) by pushing the
framings via $p$ as in Section \ref{sec:1}. The matrix of the form
$\tb_x$ in the basis formed by the connected components of $M$ is
given by:
\[[\tb_x]=\left\{\begin{array}{ccl}
\left(\begin{array}{rr}-\frac{m}{4}&-\frac{m}{4}\\-\frac{m}{4}&-\frac{m}{4}\end{array}\right)&,&
\Rr E\mbox{ 2-sided on }\Rr X'\\
(-m)&,& \Rr E\mbox{ 1-sided on }\Rr X'.
\end{array}\right.\]
\label{contract}
\end{proposition}

\begin{proof}
The point $x$ is a Q-singularity because a sufficiently small
neighbourhood $\Cc U'$ of $\Cc E$ in $\Cc X'$ has a boundary which
is a $\Qq$-homology sphere.

Let $U$ denote $p(\Cc U')$. Let us consider the map $c':\Cc
U'\longrightarrow \overline{U}$ of quotient by conjugation.

Assume first that $\Rr E$ is 2-sided on $\Rr X'$. For $i=1,2$ put
$\gamma_i=\partial\Cc U'\bigcap\Rr X'$; $\gamma_1$ and $\gamma_2$
are the two connected components of a small neighbourhood of the
cycle $\Rr E$ in $\Rr X'$. Let $A_i$ be the annulus on $\Rr X'$
bounded by $\gamma_i$ and $\Rr E$. Assume that the canonical
orientation of the disk $\overline{E}=\Cc E/\conj$ agrees with
$A_1$. Hence, since the surface $\Rr X'$ is nonorientable, the
orientation of $\overline{E}$ also agrees with the orientation of
$A_2$. Now consider the topological disks $\Sigma_1=A_1\cup_{\Rr
E}\overline{E}$ and $\Sigma_2=A_2\cup_{\Rr E}\overline{E}$ in
$\overline{U}'$.  Then it follows from Equation \ref{selfint},
Section \ref{sec:0}, that:
\[\overline{tb}_x(\gamma_1,\gamma_1)=\langle\Sigma_1,\Sigma_1\rangle_{(\overline{U}'\!,\partial\overline{U}')}
=\frac{1}{2}\langle\Cc E,\Cc E\rangle_{(\Cc U'\!,\partial\Cc U')}
=-\frac{m}{2};\]
\[\overline{tb}_x(\gamma_2,\gamma_2)=\langle\Sigma_2,\Sigma_2\rangle_{(\overline{U}'\!,\partial\overline{U}')}
=\frac{1}{2}\langle\Cc E,\Cc E\rangle_{(\Cc U'\!,\partial\Cc U')}
=-\frac{m}{2};\]
\[\overline{tb}_x(\gamma_1,\gamma_2)=\langle\Sigma_1,\Sigma_2\rangle_{(\overline{U}'\!,\partial\overline{U}')}
=\frac{1}{2}\langle\Cc E,\Cc E\rangle_{(\Cc U'\!,\partial\Cc U')}
=-\frac{m}{2}.\] First part of the proof follows from Equation
\ref{tbbar2rb}, Section \ref{sec:1}.

Now assume that the cycle $\Rr E$ is 1-sided on $\Rr X'$. In this
case, a sufficiently small neighbourhood of $\Rr E$ in $\Rr X'$ is
a M\"{o}bius band, say $A$, thus has connected boundary, say
$\gamma$. Observing that $\gamma=2\Rr E$, define $\Sigma$ to be
the relative 2-chain $A\cup_{\Rr E}2\overline{E}$ in
$\overline{U}'$. It follows that:
\[\overline{tb}_x(\gamma,\gamma)=\langle\Sigma,\Sigma\rangle_{(\overline{U}'\!,\partial\overline{U}')}
=\frac{1}{2}\cdot 4\langle\Cc E,\Cc E\rangle_{(\Cc
U'\!,\partial\Cc U')} =-2m.\]
\end{proof}

\begin{corollary}
In the previous proposition, let $\Rr X'$ be orientable. Then
\[\tb_x = \left(\begin{array}{rr}
-\frac{m}{4}&\frac{m}{4}\\\frac{m}{4}&-\frac{m}{4}
\end{array}\right).\]
\label{contract2}
\end{corollary}

\begin{proof}
We should only note that in the case of orientable $\Rr X'$, $\Rr
E$ is 2-sided on $\Rr X'$ and the orientation of $\overline{E}$
does not agree with the orientations of $A_1$ and $A_2$
simultaneously.
\end{proof}

\begin{theorem} \label{calcul} Let $x\in\Cc X$ be a normal, numerically
Gorenstein, unibranch real Q-singularity;
$\Cc U$ a small  compact cone-like real neighbourhood
of $x$; $\Rr M=\Rr X\cap\partial\Cc U$; $\Gamma$ be the graph
of a very good resolution $\rho:\Cc\widetilde{U}\rightarrow
\Cc U$ at $x$;
$\Rr\widetilde{U}$ be the preimage $\rho^{-1}(\Rr U)$;
${\mathcal W}_{\Rr}$ be the characteristic set on $\Gamma_{\Rr}$ and
$N$ be the number of vertices on $\Gamma_{\Rr}$. Then:
\begin{equation} \tb(\Rr M)=N-1
+\sum_{\scriptstyle 1\leq i\leq s \atop
\scriptstyle e_i\in {\mathcal W}_{\Rr}} n'_i.
\label{eqCalcul}\end{equation}
\end{theorem}

\begin{proof}
First note that the normality of the singularity $x$ implies that
$\Gamma$ is connected.
Since $x$ is assumed to be numerically Gorenstein,
the canonical class of $\Cc\widetilde{U}$ is integral so that one can define the characteristic set
of the graph. Furthermore since $x$ is a Q-singularity, the graph $\Gamma$
is a tree and each exceptional curve is a sphere (see e.g. \cite{Dim}).

The surface $\Rr\widetilde{U}$ is smooth, connected and with
boundary $\Rr M$. $\Rr\widetilde{U}$ need not to be orientable. It
is nonorientable if and only if $n_i$ ($1\leq i\leq s$) is odd for
at least one real exceptional sphere (see \cite{Sil}, Proposition
II.6.9 and Corollary II.6.10). If $n_i$ is even for all
exceptional spheres, then $\Rr\widetilde{U}$ is oriented and
$c_1(\Cc\widetilde{U})\Big|_{\Rr\widetilde{U}}$ is even, i.e. by
definition $\mathcal W_{\Rr}$ is empty, so that the second term on
the right hand side of Equation \ref{eqCalcul} disappears and one
obtains Equation \ref{selfint}, Section \ref{sec:0}:
\[\tb(\Rr M)=-\chi(\Rr\widetilde{U}) = N-1.\]

In the general case where $\Rr\widetilde{U}$ is nonorientable,
we construct a real algebraic surface $\Cc\widehat{U}$ with an oriented
real part by getting rid of the terms that appear in
$c_1(\Cc\widetilde{U})\Big|_{\Rr\widetilde{U}}$.
So, let \[F=\{E_i|e_i\in{\mathcal W}_{\Rr}\}\bigcup \{E_i|e_i\in{\mathcal C}_j
\mbox{ for some }j\,(1\!\leq\!j\!\leq\!s), \mbox{ such that }
e_j\in{\mathcal W}_{\Rr}\}.\]
We define: \[\pr:\Cc\widetilde{U}\rightarrow\Cc\widehat{U}=
\Cc\widetilde{U}\big/F.\]

It is well-known that $c_1(\Cc\widehat{U})|_{\Rr\widehat{U}}$
is even if and only if $w_1(\Rr\widehat{U})$ is zero.
Since each exceptional curve in the canonical class of $\Cc U$
is contracted in $\Cc\widehat{U}$, the real part $\Rr\widehat{U}$
is orientable. However, $\Rr\widehat{U}$ is not smooth.
Let the number of singularities of $\Cc\widehat{U}$
be $t$ which is exactly the number of vertices in ${\mathcal W}_{\Rr}$. Each
singular point is real and is the image $\pr(E_i)$ where $E_i$ corresponds to the vertex
$e_i\in{\mathcal W}_{\Rr}$; we denote this singular point by $x_i$.

\begin{lemma} Each singular point of $\Cc\widehat{U}$ is a Q-singularity.
\end{lemma}

\begin{proof}
Let $x$ be such a singularity.
The subgraph of $\Gamma$ spanned
by the vertices corresponding to the exceptional spheres in the
preimage $\pr^{-1}(x)$
is a connected, negative definite tree, i.e. a small regular
neighbourhood of $\pr^{-1}(x)$ in $\Cc\widetilde{U}$ has boundary
a $\Qq$-homology sphere.
\end{proof}

\begin{lemma} Let $x_i\in\Rr\widehat{U}$ be one of above singularities. In the basis of
the homology classes of the real branches of $x_i$,
\[[\tb_{x_i}]= \left\{\begin{array}{ccl}
\left(\begin{array}{rr}\frac{n'_i}{4}&\frac{n'_i}{4}
\\\frac{n'_i}{4}&\frac{n'_i}{4}\end{array}\right)&,&
\Rr E\mbox{ 2-sided on }\Rr X'\\
(n'_i)&,& \Rr E\mbox{ 1-sided on }\Rr X'.
\end{array}\right.\]
\label{lemContr}
\end{lemma}

\begin{proof} Let $\Cc\widehat{U}_{x_i}\subset\Cc\widehat{U}$ be a
sufficiently small compact cone-like neighbourhood of $x_i$ and
$\Cc M_{x_i}$ be the boundary $\partial\Cc\widehat{U}_{x_i}$. Let
${\mathcal E}$ denote the subspace of
$H_2(\pr^{-1}(\Cc\widehat{U}_{x_i}))$ generated by the homology
classes of exceptional spheres in $\pr^{-1}(x_i)$. The complex
intersection form in $\Cc\widehat{U}_{x_i}$ is the restriction of
the intersection form in $\pr^{-1}(\Cc\widehat{U}_{x_i})$ to the
orthogonal complement ${\mathcal E}^{\perp}$ of ${\mathcal E}$ in
$H_2(\pr^{-1}(\Cc \widehat{U}_{x_i}))$, since the restriction of
the intersection form in $\pr^{-1}(\Cc\widehat{U}_{x_i})$ to
${\mathcal E}$ is nondegenerate (see \cite{Fin2}, Section 5.2).
The result of the corresponding orthogonalisation process on the
matrix of the intersection form coincides with the continued
fractions of the self intersection numbers of the vertices on the
contracted arms. The self intersection of $[E_i]$ restricted to
the subspace $\langle {\mathcal E}^{\perp},[E_i]\rangle$ is
therefore equal to the number $n'_i$ defined in Section
\ref{sec:2}. Since the singularity $x_i$ is obtained by
contracting $E_i$, the result follows from Proposition
\ref{contract}.
\end{proof}

Now, since $\Rr\widehat{U}$ is oriented, by Equation
\ref{selfint}, Section \ref{sec:0}, we get:
\[\begin{array}{ll}
\tb(\Rr M)=\tb(\Rr\widehat{M})&
=\langle\Rr\widehat{U},\Rr\widehat{U}\rangle_{(\Cc\widehat{U}\!,\partial\Cc\widehat{U})}\\
&=-\chi\big(\Rr\widehat{U}-\{t \mbox{ points}\}\big)
+{\displaystyle \sum_{\scriptstyle 1\leq i\leq s \atop
\scriptstyle e_i\in {\mathcal W}_{\Rr}}}
\langle\Sigma_i,\Sigma_i\rangle_{(\Cc\widehat{U}_{x_i}\!,\partial\Cc\widehat{U}_{x_i})}
\end{array}\] where $\Sigma_i$ is the connected piece of the
surface $\Rr\widehat{U}$ bounded by the connected components of
the real link $\Rr\widehat{M}_{x_i}$ of ${x_i}$ (number of
connected components is either one or two).


For any $e_i\in{\mathcal W}_{\Rr}$, first let the curve $\Rr E_i$ be 2-sided
on the nonorientable surface  $\Rr\widetilde{U}$ and $\gamma_{i1}$, $\gamma_{i2}$
be the two connected components of $\Rr\widehat{M}_{x_i}$. It easy to see that
\[\begin{array}{ll}
\langle\Sigma_i,\Sigma_i\rangle_{(\Cc\widehat{U}_{x_i}\!,\partial\Cc\widehat{U}_{x_i})}
&=\tb_{x_i}(\gamma_{i1},\gamma_{i1})+\tb_{x_i}(\gamma_{i2},\gamma_{i2})+2\tb_{x_i}(\gamma_{i1},\gamma_{i2})\\
&=4\cdot\frac{n'_i}{4}=n'_i \end{array}\] where the second
equality follows from Lemma \ref{lemContr}. Similarly, letting
$\Rr E_i$ be 1-sided on $\Rr\widetilde{U}$, we get ${\displaystyle
\langle\Sigma_i,\Sigma_i\rangle_{(\Cc\widehat{U}_{x_i}\!,\partial\Cc\widehat{U}_{x_i})}=n'_i.}$
Hence,
\[\begin{array}{ll}\tb(\Rr M)&=-\chi\big(\Rr\widehat{U}-\{t \mbox{ points}\}\big)
+{\displaystyle \sum_{\scriptstyle 1\leq i\leq s \atop
\scriptstyle e_i\in {\mathcal W}_{\Rr}}} n'_i\\
&=-(2-N-1)
+{\displaystyle \sum_{\scriptstyle 1\leq i\leq s \atop
\scriptstyle e_i\in {\mathcal W}_{\Rr}}} n'_i\\
&=N-1
+{\displaystyle \sum_{\scriptstyle 1\leq i\leq s \atop
\scriptstyle e_i\in {\mathcal W}_{\Rr}}} n'_i.\end{array}\]

\vspace{-0.3cm}

\end{proof}
Note that the above calculation is an application of the Integral
Formula (see \cite{Fin1}) applied to dimension two.

\begin{figure}[h]
\begin{center}

\hspace{1.2cm} $-2$ \hspace{.4cm} $-2$ \hspace{.5cm} $-2$
\hspace{.5cm} $-2$ \hspace{.6cm} $-3$

\vspace{.2cm}

$-3$ \hspace{.18cm} $-2$ \hspace{5.2cm}

\vspace{1.3cm}

\hspace{1.2cm} $-2$ \hspace{.4cm} $-2$ \hspace{.5cm} $-2$
\hspace{.5cm} $-2$ \hspace{.6cm} $-3$

\vspace{-2.3cm}

\epsfig{file=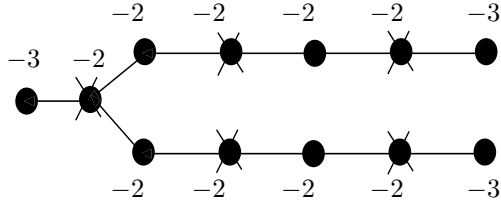,height=2cm,width=6.5cm}

\vspace{.5cm}

\caption{Resolution graph for $x^5+y^8+z^2=0$} \label{gra58}
\end{center}
\end{figure}

\noindent \textbf{Example:} Consider the singularity on the
surface $\Cc X_{\pm}=\{x^5+y^8\pm z^2=0\}$ with the resolution
graph $\Gamma(5,8)$ as in Figure \ref{gra58}. The five vertices
belonging to the characteristic set ${\mathcal W}$ are marked with
a cross in the figure.

The complex conjugation $\conj_-$ on the resolved surface
$\Cc\widetilde{X}_-$ fixes all vertices of $\Gamma(5,8)$
(Proposition \ref{prop0}). Therefore ${\mathcal W_{\Rr}}$ is the
set of the five vertices marked. Meanwhile, for
$\Cc\widetilde{X}_+$, ${\mathcal W_{\Rr}}$ contains only the
rupture vertex. In both cases $\Rr\widetilde{X}_{\pm}$ is
nonorientable.

It follows now from Theorem \ref{calcul} that
\[\tb_-=12-1+5\cdot(-2)=1 \mbox{ and }
\tb_+=2-1+\big(-2-2\cdot(-\frac{9}{11})\big)=\frac{7}{11}.\]
The values for
$\tb_{\pm}$ can be calculated by the method in \cite{Fin2}
as well giving the same results.

\subsection{Computation of the linking form for singularities with
more than one real branch: an example} \label{sec:7}
Suppose now that the singularity $x\in\Cc X$ has $\rho$ real
branches with $\rho>1$. Then the real part $\Rr M$ of the link of
singularity is disconnected. Hence the linking form $\tb$ is given
by $\rho\times\rho$ matrices. Let us take for example the complex
surface $\Cc X$ defined as the double cover of $\Cc^2$ branched
along $\{y(x^5+y^4)=0\}$. The boundary of the surface $\Rr U$ has
two components $\alpha$ and $\beta$. With the methods of
\cite{Fin2}, one can compute that the $2\times 2$ matrix $[\tb]$
corresponding to $\tb$ in the basis $\{[\alpha],[\beta]\}$ is:

\[[\tb]=\left(\begin{array}{rr}{\displaystyle \frac{3}{2}}
&{\displaystyle -\frac{5}{2}}
\\{\displaystyle -\frac{5}{2}}&
{\displaystyle \frac{3}{2}}\end{array}\right)\]

\begin{figure}[h]
\begin{center}
$\hspace{.1cm} -2 \hspace{.3cm} -2\hspace{.4cm} -2\hspace{.3cm}
-2\hspace{.4cm} -2\hspace{.3cm} -2\hspace{.3cm} -2\hspace{.4cm}
-2\hspace{.3cm} \;-2\hspace{.4cm} -3$
\epsfig{file=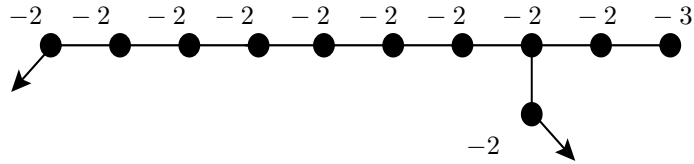,height=1.8cm,width=9cm} \caption{Resolution
graph for the singularity $z^2=y(x^5+y^4)$} \label{z2yx5y4}
\vspace{-1.4cm}

\hspace{3.5cm} $-2$

\vspace{1.2cm}

\end{center}
\end{figure}

Alternatively, one can compute this matrix using the arguments in
the proof of Theorem \ref{calcul}. The very good resolution graph
$\Gamma$ corresponding to the singularity at $0$ is as in Figure
\ref{z2yx5y4} where the arrows stand for the two components of the
proper transform of the branching locus.

\begin{figure}[h]
\begin{center}
\epsfig{file=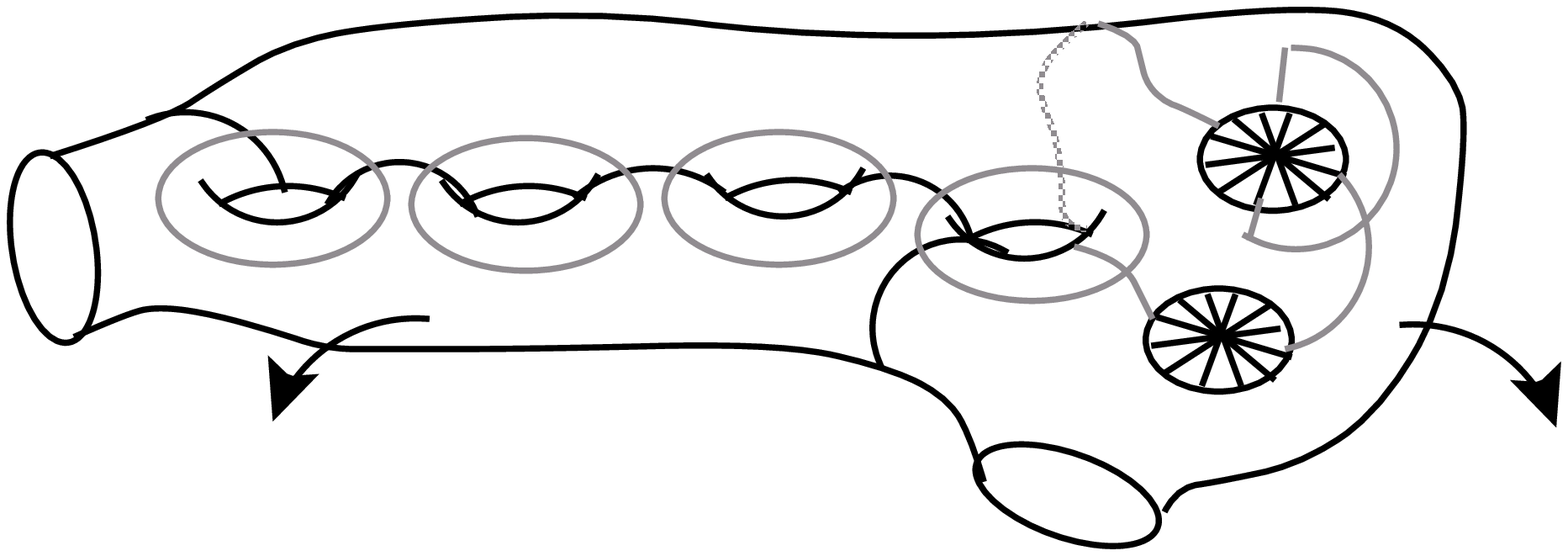,height=3cm,width=7cm} \caption{Real part of
the surface $\{z^2=y(x^5+y^4)\}$ after resolution}
\label{surz2yx5y4}

\vspace{-3.9cm}

\hspace{3.cm}$\gamma_6$

\hspace{5.3cm}$\gamma_7$


\vspace{-.4cm} $\gamma_1$ \hspace{5.3cm}

$\gamma_2\;\;\;\;\;\;\;\;\;\gamma_3\;\;\;\;\;\;\;\;\;\gamma_4$\hspace{1.cm}

\vspace{.6cm}

$\alpha$\hspace{3.6cm}$\gamma_5$\hspace{3.4cm}

\vspace{.5cm}

\hspace{1cm}$\Sigma_1$\hspace{5.2cm}$\Sigma_2$

\hspace{1.5cm} $\beta$


\end{center}
\end{figure}

There may be assigned two real structures on $\Cc\widetilde{X}$.
Each of the complex conjugations fixes each vertex on the
resolution graph. Fix any one of the complex conjugations. Then
the surface $\Rr\widetilde{U}$ is a non-oriented surface with the
two boundary components $\alpha$ and $\beta$ (see Figure
\ref{surz2yx5y4}).

Note that the curves $\gamma_1,\ldots,\gamma_5$ separate
$\Rr\widetilde{U}$ into two disjoint surfaces; let us denote by
$\Sigma_1$ the one which contains $\alpha$ in its boundary and by
$\Sigma_2$ the other surface. Note that $\Sigma_1$ is orientable
but $\Sigma_2$ is not.

The curves are the real parts of the $-2$-spheres
$E_1,\ldots,E_5$. Let ${\mathcal S}={e_1,\ldots,e_5}$. Furthermore
${\mathcal W}=\{e_5,e_6\}$ is a characteristic set of $\Gamma$.

Similarly as in the proof of Theorem \ref{calcul}, let
$F=\{E_1,\ldots,E_5,E_6\}$ and let
$\pr:\Cc\widetilde{U}\rightarrow\Cc\widehat{U}=\Cc\widetilde{U}\big/F$.
Let us denote the Q-singularities $\pr(E_i)\in\Cc\widehat{U}$
$(i=1,\ldots 6)$ by $x_i$. The surface
$\widehat{\Sigma}_1=\pr(\Sigma_1)$ is orientable, has boundary
$\pr(\alpha)$ and contains five Q-singularities $x_1,\ldots,x_5$.
The surface $\widehat{\Sigma}_2=\pr(\Sigma_2)$ is orientable, has
boundary $\pr(\beta)$ and contains six Q-singularities
$x_1,\ldots,x_5,x_6$. It follows by the proof of Theorem
\ref{calcul} and Proposition \ref{contract} that
\[\begin{array}{ll} \tb(\alpha,\alpha)
=\langle\widehat{\Sigma}_1,\widehat{\Sigma}_1\rangle_{\Cc\widehat{U}\!,\partial\Cc\widehat{U}}
&\displaystyle =-\chi(\widehat{\Sigma}_1-\{x_1,\ldots,x_5\})+5\cdot(-\frac{1}{2})=4-\frac{5}{2}=\frac{3}{2};\\
\tb(\beta,\beta)
=\langle\widehat{\Sigma}_2,\widehat{\Sigma}_2\rangle_{\Cc\widehat{U}\!,\partial\Cc\widehat{U}}
&\displaystyle =-\chi(\widehat{\Sigma}_2-\{x_1,\ldots,x_6\})+5\cdot(-\frac{1}{2})+(-2)\\
&\displaystyle =6-\frac{5}{2}-2=\frac{3}{2};\\
\tb(\alpha,\beta)=
\langle\widehat{\Sigma}_1,\widehat{\Sigma}_2\rangle_{\Cc\widehat{U}\!,\partial\Cc\widehat{U}}
&\displaystyle
=5\cdot(\frac{-1}{2})=-\frac{5}{2}=\tb(\beta,\alpha).\end{array}\]

\section{On Thurston-Bennequin numbers of Brieskorn double points}
\label{sec:8}
The Thurston-Bennequin numbers corresponding to the Brieskorn
double points $x^m+y^n\pm z^2=0$ ($\gcd(m,n)=1$) will be denoted
by $\tb_\pm(m,n)$.

\begin{corollary} \label{q-int}
The number $\tb_-(m,n)$ is always an integer.
\end{corollary}
\begin{proof} By Proposition \ref{prop0}, Section
\ref{sec:4},  all vertices of $\Gamma(m,n)$ are fixed under the
$\conj_-$, i.e. $\Gamma_{\Rr}=\Gamma$, so that all the non-integer
terms in Equation \ref{eqCalcul}, Section \ref{sec:6}, disappear.
\end{proof}

The Thurston-Bennequin number $\tb_+ (m,n)$ is in
general non-integral.

\begin{proposition}
\label{e0notinW}
Let $m$ be even and $n$ odd. The rupture vertex $e^0$ is
in ${\mathcal W}$ if and only if $4\!\not|m$.
\end{proposition}

\begin{proof}
The vertex $e_0=conj_{\pm}(e^0)$ in the graph $\Gamma_{f}$ of
$\Cc\widetilde{Y}$ has even multiplicity $mn$. Therefore by
Proposition \ref{b0cift} and Proposition \ref{tekcift}, Section
\ref{sec:5},  the proof follows.
\end{proof}

\begin{corollary}
\label{q+int} Assume that $4|m$ and $n$ is odd. Then
$\tb_+(m,n)$ is integer.
\end{corollary}

\begin{proof}
From the proposition above, it follows that
there is no (non-integer) contribution of
imaginary arms.
\end{proof}

The theorem below shows that if the value of $tb_{\pm}(m,n)$ is known for
sufficiently small values of $m$ and $n$, the number $tb_{\pm}(m',n')$ can be
computed for greater integers $m'$ and $n'$.

\begin{theorem}
\label{period}  If $m$ is odd then \[\tb_{\pm}(m, n+4m)=\tb_{\pm}(m,n).\]
If $4|m$ then \[\tb_{\pm}(m, n+2m)=\tb_{\pm}(m,n).\]
If $m\equiv 2\mod\,4$ then
\[\tb_-(m, n+2m)-\tb_-(m,n)=4,\]
\[\tb_+(m, n+2m)-\tb_+(m,n)=\frac{4}{n(n+2m)}.\]
\end{theorem}

\begin{proof}
Assume first that $m$ is odd. By Corollary \ref{cor1}, Section
\ref{sec:4},  there are two vertices appended to each ($n$)-arm of
$\Gamma(m,n)$ to construct $\Gamma(m,n+4m/(m,2))$; let us call
them $e_1$ and $e_2$ with $e_1$ being the outer-most vertex. Take
an ($n$)-arm of $\Gamma(m,n)$; let us call the two vertices to be
appended as $e_1$ and $e_2$ with $e_1$ being the outer-most
vertex. By Proposition \ref{prop2}, Section \ref{sec:4}, $e_1$ has
self intersection $-2$.

To prove the theorem, it is enough to show that $e_1\in{\mathcal
W}$ and $e_2\not\in{\mathcal W}$ since in that case it follows
from Theorem \ref{calcul}, Section \ref{sec:6},  that
\[\tb_-(m, n+4m)=\tb_-(m,n)-(-2)+(-2)=\tb_-(m,n).\]

Recall that by Proposition \ref{prop2}, the terminal part  added
to $\Gamma_f$ to construct the graph $\Gamma_{f'}$ corresponding
to $\widetilde{Y}'$ is as in Figure \ref{4nodes}(a). The
exceptional curves $E_1'$ and $E_2'$ have odd multiplicities $m$
and $3m$ respectively. Therefore by Proposition \ref{lem3},
Section \ref{sec:5},  the vertices of $\Gamma(m,n+4m)$
corresponding to the exceptional curves $\widetilde{p}^{-1}(E'_1)$
and $\widetilde{p}^{-1}(E'_2)$ are not in ${\mathcal W}$. It
follows that the vertex corresponding to the exceptional curve
$\widetilde{p}^{-1}(E_1)$ must appear in ${\mathcal W}$ for
Equation \ref{defw}, Section \ref{sec:5}, to be satisfied in
$\Cc\widetilde{U}'$. Again in order to satisfy Equation
\ref{defw}, the vertex corresponding to the exceptional curve
$E_2$ is not included in ${\mathcal W}$.

Assume now that $m$ is divisible by 4 and $n$ is odd. By
Proposition \ref{prop0}, Corollary \ref{cor1}, Section
\ref{sec:4}, and Proposition \ref{e0notinW}, it follows that
$\tb_+(m, n+2m)=\tb_+(m,n)$. For the relation on the invariant
$\tb_-$, we first note that, as in the proof of Proposition
\ref{prop1}, Section \ref{sec:4},  the terminal part of the
resolution graph $\Gamma_{f'}$  is as shown in Figure \ref{2node}.
The multiplicity of the sphere $E_1$ is $m$ hence by Lemma
\ref{mel}, the preimage $(\widetilde{p})^*(E_1)$ does not
contribute mod $2$ to
$(\widetilde{p})^*(\frac{1}{2}[\widehat{A'}])$.

\begin{figure}[h]
\begin{center}

$e_1\;\;\;\;\;e_2$ \hspace{1.2cm}

\epsfig{file=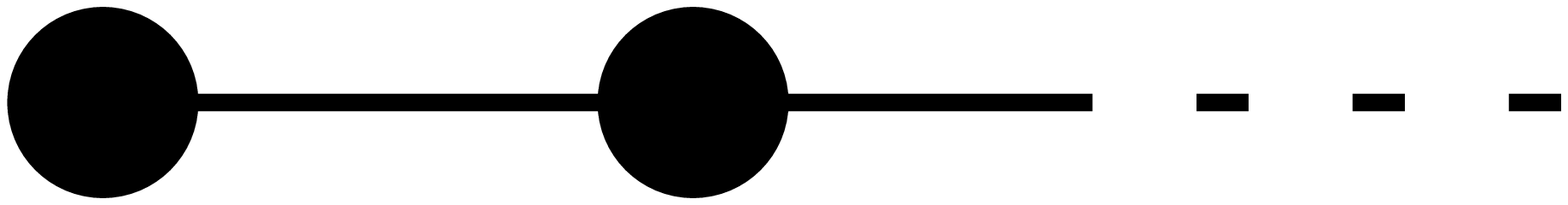,height=0.35cm,width=2.8cm}

${\displaystyle {m\atop -2}\;\;\;\,{2m\atop-k}}$ \hspace{1.2cm}

\caption{}\label{2node}

\end{center}
\end{figure}

The coefficient $b_1$ of the cycle $E_1$ in
$c_1(\Cc\widetilde{V}')$ is $-1$ because $E_1$ is produced after
the first blow-up of $\Cc^2$ (cf. Proposition \ref{lem4}, Section
\ref{sec:5}). Therefore by Proposition \ref{mel}, Section
\ref{sec:5}, the vertex of $\Gamma(m,n+4m)$ corresponding to
$\widetilde{p}^{-1}(E_1)$ is in ${\mathcal W}$. As a consequence,
$\widetilde{p}^{-1}(E_2)$ is not in ${\mathcal W}$.

(One can also show that $\widetilde{p}^{-1}(E_2)$ is not in ${\mathcal W}$
by the following reasoning:
the multiplicity of $E_2$ is $2m$. Since $E_2$ is produced after
the second blow-up of $\Cc^2$, the coefficient $b_2$ is $-2$ so
that $a_2=0 \mbox{ mod }2$. Therefore the vertex of $\Gamma(m,n+4m)$
corresponding to  $\widetilde{p}^{-1}(E_2)$ is
not in ${\mathcal W}$.)

Now assume that $m\equiv 2$ $(\!\!\!\mod\,4)$ and $n$ is odd.
The terminal
part of the resolution graph $\Gamma_{f'}$  is again as in Figure
\ref{2node}. Note that $b_1=-1$ and $b_2=-2$ so that by Proposition \ref{tekcift},
$e_1,e_2\not\in{\mathcal W}$. It follows from Theorem \ref{calcul} that
\[tb_-(m, n+2m)-\tb_-(m,n)=-\chi(\Rr X')+\chi(\Rr X)=2\cdot 2=4.\]

For the last claim in the theorem, we observe by Proposition
\ref{e0notinW} that $e^0\in{\mathcal W}_{\Rr}$. It follows from
Corollary \ref{cor1} that $\Gamma(m,n+2m)$ and $\Gamma(m,n)$
differ in the pair of $(n+2m)$-arms of the former and pair of
($n$)-arms of the latter. Let $\sigma'$ be an $(n+2m)$-arm of
$\Gamma(m,n+2m)$ and $\sigma$ be an ($n$)-arm of $\Gamma(m,n)$. By
Theorem \ref{calcul}
\[tb_+(m, n+2m)-\tb_+(m,n)=2(\frac{1}{n^{\sigma'}}-\frac{1}{n^{\sigma}}).\]
By \cite{OrW} Theorem 3.6.1, it follows that
\[\frac{1}{n^{\sigma'}}-\frac{1}{n^{\sigma}}=-\frac{1}{m(n+2m)}+\frac{1}{mn}=\frac{2}{n(n+2m)}\]
and the proof follows.
\end{proof}

\begin{theorem}
\label{simetri} Let $m$ be odd and $k,t\in\Zz^+$, $4km>t$. Then
\[\tb_{\pm}(m,4km-t)+\tb_{\pm}(m,t)=-2.\]
\end{theorem}

\begin{proof}
Consider the weighted
homogeneous surface  $\Cc P(4k,1,1)$ with a unique singularity at
$[1:0:0]$ (cf. \cite{Dim} Appendix B for a
brief introduction to weighted projective spaces)
and the curve $\Cc A=\{x^m+u^{4km-t}v^{t}=0\}$ in $\Cc P(4k,1,1)$.
$\Cc A$ has a pair of
singularities at $[0:0:1]$ and $[0:1:0]$. Let
$\Cc\widetilde{P}(4k,1,1)$ obtained by blowing-up $\Cc P(4k,1,1)$
at $[1:0:0]$.  Denote the resolution by $\rho$.

The plan of the proof is as the following. First we construct
a  branched double covering $p:\Cc X\rightarrow\Cc P(4k,1,1)$
whose branching locus contains the curve $\Cc A$. The surface $\Cc X$
has the real Brieskorn double points corresponding to $(m,t,2)$ and $(m,4km-t,2)$.
Then we show that $\Rr X$ is orientable and $\overline{X}=\Cc X/\conj$ is a $\Qq$-homology sphere
so that $\big\langle \Rr X,\Rr X\big\rangle_{\overline{X}}=0$.
Using this identity, we obtain a relation between $\tb_{\pm}(m,4km-t)$ and $\tb_{\pm}(m,t)$.

The smooth surface
$\Cc\widetilde{P}(4k,1,1)$ is biholomorphic to the
Hirzebruch surface $\Ff_{4k}$  and
the section $S_\infty\subset\Cc\widetilde{P}(4k,1,1)$
has self intersection $-4k$ (see e.g. \cite{Bea}). Let $S_0$ be the zero section  and $F$ be
a fiber of ruling in $\Cc\widetilde{P}(4k,1,1)$. One can show that
$[\rho^{-1}(\Cc A)]\!=\!m[S_0]$ and $\big([\rho^{-1}(\Cc A)]+[S_\infty]\big)\!=\!(m+1)[S_0]-4k[F]$.
Hence, $[\rho^{-1}(\Cc A)]$ is odd and $\big([\rho^{-1}(\Cc A)]+[S_\infty]\big)$ is even.
Therefore there exist a double covering
$\widetilde{p}:\Cc\widetilde{X}\rightarrow\Cc\widetilde{P}$
branched along $\rho^{-1}(\Cc A)\cup S_\infty$. The covering
$\widetilde{p}$ is associated
either with the complex conjugation $\conj_-$ or with $conj_+$.
Let us assume that $\widetilde{p}$ is associated with $\conj_-$ first.

Consider the commutative diagram:
\[
\begin{array}{rcccl}
\Cc\widetilde{X}&\stackrel{\textstyle\widetilde{\rho}}{\longrightarrow}&\Cc
X&\stackrel{\textstyle c'}{\longrightarrow}&\Cc
X/\conj=\overline{X}\\ \widetilde{p}\Big\downarrow&&\Big\downarrow
p&&\Big\downarrow\overline{p}\\ \Ff_{4k}=
\Cc\widetilde{P}(4k,1,1)&\stackrel{\textstyle\rho}{\longrightarrow}&\Cc
P(4k,1,1)&\stackrel{\textstyle c}{\longrightarrow}&\Cc
P(4k,1,1)/\conj=\overline{P}
\end{array}
\]
where $p:\Cc X\rightarrow\Cc P(4k,1,1)$ is the double covering with
branching locus $\Cc A\cup\rho(S_\infty)=\Cc A\cup\{[1:0:0]\}\subset\Cc P(4k,1,1)$
and $c$, $c'$ are maps of quotient by conjugation.
Note that $\Cc X=\Cc\widetilde{X}/E$ where $E=\widetilde{p}^{-1}(S_\infty)$
with $E^2=-2$
and $\Cc X$ has isolated singularities at $p^{-1}([1:0:0])$,
$p^{-1}([0:1:0])$ and $p^{-1}([0:0:1])$,
the last two being Brieskorn
double points corresponding to $(m,t,2)$ and $(m,4km-t,2)$ respectively.

It is easy to check that $c_1(\Cc\widetilde{X})$ is even. Therefore
the surface $\Rr\widetilde{X}$ is orientable.
In particular, $\Rr\widetilde{X}$ is a torus and $\Rr X$ is a pinched torus.
To see this, let $\Delta$ be the subset of $\Rr\widetilde{P}(4k,1,1)$
such that $\Delta=\widetilde{p}(\Rr\widetilde{X})$. Then
$\Delta$ is the closure of one of the connected components of
$\Rr\widetilde{P}(4k,1,1)-\big(\rho^{-1}(\Rr A)\cup S_\infty\big)$.
Therefore $\Delta$ is an annulus (see Figure \ref{silindir}).
Note that the choice of $\Delta$
is determined by the choice of the complex conjugation $\conj_{\pm}$.

\begin{figure}[h]
\begin{center}
\epsfig{file=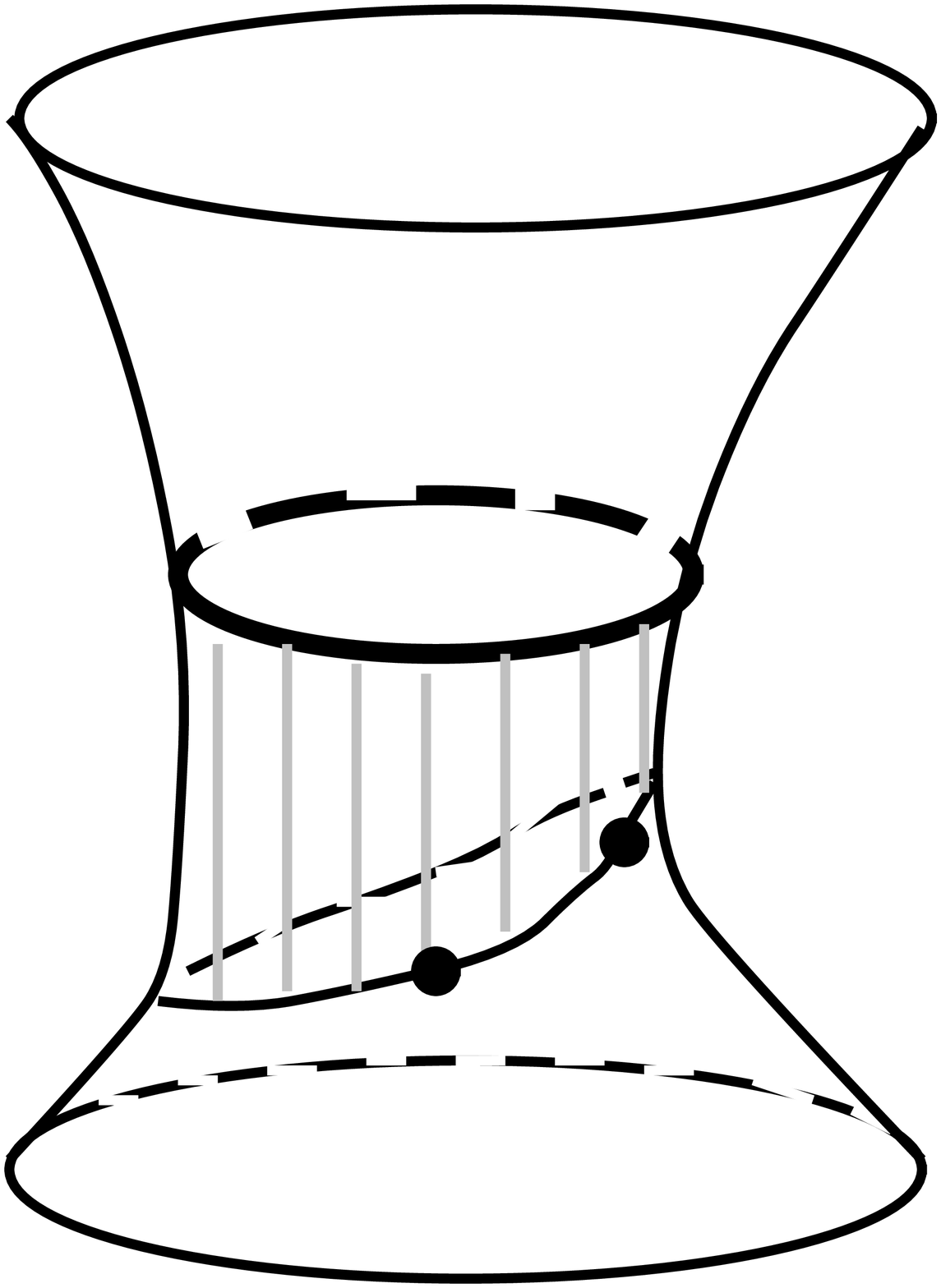,height=4cm,width=3cm}
\caption{$\Rr\widetilde{A}\cup\Rr S_\infty$ bounds an annulus on $\Rr
\widetilde{P}(4,1,1)$} \label{silindir}

\vspace{-2.9cm}

\hspace{2.4cm} $\Rr S_\infty$

$\Delta$

\vspace{0.3cm}

$\rho^{-1}(\Rr A)$ \hspace{3.5cm}

\vspace{1cm}

\end{center}
\end{figure}

$\Rr X$ has three singular points: $P_1=p^{-1}([0:1:0])$,
$P_2=p^{-1}([0:0:1])$ and $P_3=p^{-1}([1:0:0])$. The point $P_3$
is obtained by contracting $E$ in $\Cc\widetilde{X}$. It follows
from Corollary \ref{contract2}, Section \ref{sec:6}, that

\begin{equation}
[\tb_{P_3}]=
\left[ \!\!\begin{array}{rr}
-\frac{k}{2}&\frac{k}{2}\\\frac{k}{2}&-\frac{k}{2}\end{array}\right].
\label{konematris} \end{equation}

Finally we observe that $\overline{P}$ and $\overline{X}$ are $\Qq$-homology 4-spheres.
For the former, we note the well-known fact that
$\Cc P(k,1,1)$ is a $\Qq$-homology manifold
with cohomology
$H^*(\Cc P(k,1,1);\Qq)=H^*(\Cc P^2;\Zz)\otimes\Qq$.
Then it follows
from \cite{Bre}, Section III.5, that $\overline{P}(k,1,1)$ is
a $\Qq$-homology sphere.

In order to prove that $\overline{X}$ is a $\Qq$-homology sphere,
it is enough to show that
$\overline{p}:\overline{X}\rightarrow \overline{P}$ is branched
along a 2-sphere. In fact,
the branching locus  is the {\it Arnold surface}
$\mathcal{A}=\overline{A}\cup_{\partial} c\big(\Rr P(4k,1,1)
-\mbox{int}(p(\Rr X))\big)
\subset\overline{P}(k,1,1)$ where $\overline{A}=c(\Cc A)$.
Since $\Cc{A}$ is a sphere, it follows immediately that $\mathcal{A}$
is a sphere.

Therefore there is no rational 2-homology in
$\overline{X}$. In particular, the pinched torus $c'(\Rr X)$ has
self intersection 0 in $\overline{X}$.
It follows from this fact and Equation
\ref{konematris} that:
\[
\begin{array}{ll}
0&=\big\langle c'(\Rr X),c'(\Rr X)\big\rangle_{\overline{X}}\\
&=\overline{\tb}_{P_1}(\Rr X,\Rr X)+\overline{\tb}_{P_2}(\Rr X,\Rr X)
+\overline{\tb}_{P_3}(\Rr X,\Rr X)-2\chi(\Rr X-3\mbox{ points})\\
& =\overline{\tb}_{P_1} +\overline{\tb}_{P_2}
+\big(-\frac{k}{2}+\frac{k}{2}-\frac{k}{2}+\frac{k}{2}\big)-2\cdot(-2)\\
&=\overline{\tb}_{P_1} +\overline{\tb}_{P_2}+4.
\end{array}\] The proof follows from Equation \ref{tbbar2rb}, Section \ref{sec:3}.
For an extended discussion of this proof, see also \cite{Ozt}.
\end{proof}

\begin{remark}
Theorem \ref{simetri} is not true for the two singularities
$x^m+y^{2km-t}+z^2=0$ and $x^m+y^t+z^2=0$ for arbitrary $k$ with
$k,t\in\Zz^+$. In fact, the proof fails for odd $k$ since in that
case $\widetilde{p}^{-1}(S_\infty)$ would have self intersection
$\frac{-2k}{2}=-k$ in $\Cc\widetilde{X}$, that is,
$\Rr\widetilde{X}$ and $\Rr X$ would be nonorientable so that the
Equation \ref{selfint} would not be valid.
\end{remark}

\begin{theorem}
\label{simetreven} Let $m$ be even and $k,t\in\Zz^+$, $2km>t$. Then
\[\tb_-(m,2km-t)+\tb_-(m,t)=
\left\{\begin{array}{ll}
-4, & 4|m\\
-4+4k, & m\equiv 2 \mod 4.
\end{array}
\right.\]
\end{theorem}

\begin{proof}
The idea and constructions are much similar to the previous proof.
This time let us consider $\Cc A=\{x^m+u^{2km-t}v^{t}=0\}$ in the
weighted homogeneous surface $\Cc P(2k,1,1)$. The surface $\Rr
P(2k,1,1)$ and $\Rr A$ are now as in Figure \ref{kone2}. Let us
denote by $\Delta$ the disc bounded by $\Rr A$ in $\Rr P(2k,1,1)$.
The curve $\Cc A$ is even in $H_2(\Cc P(2k,1,1);\Zz)$ and there
exists a double covering of $\Cc P(2k,1,1)$ with branching locus
$\Cc A$.

\begin{figure}[h]
\begin{center}
\epsfig{file=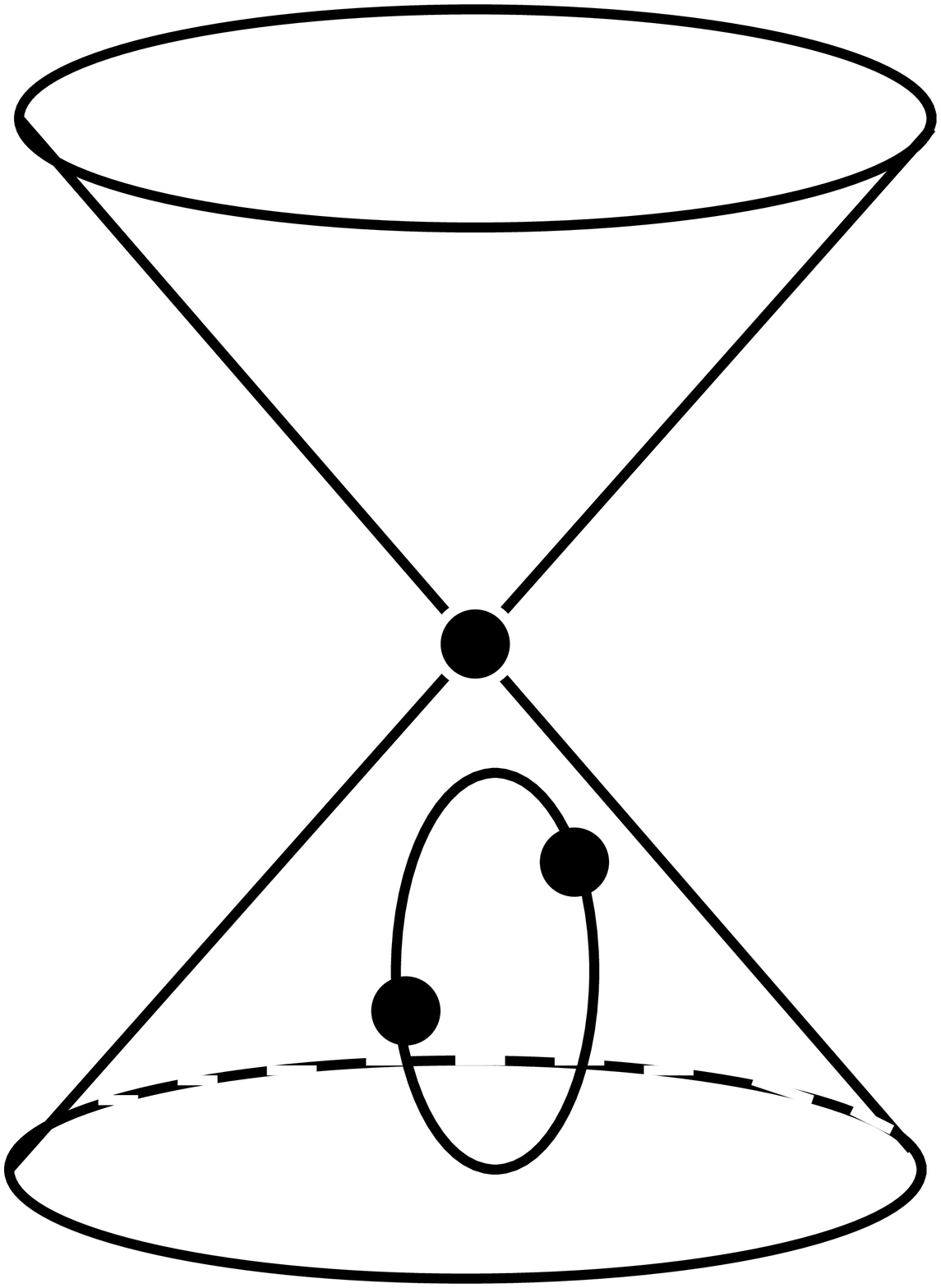,height=4cm,width=3cm} \caption{$\Rr A$
bounds a disc on $\Rr P(2k,1,1)$} \label{kone2}

\end{center}
\end{figure}

Let us first consider the case when the double covering is equipped with
the complex conjugation  $\conj_-$; denote the covering surface by $\Cc X_-$.
Then $p(\Rr X_-)=\Rr P(2k,1,1)-\mbox{int}(\Delta)$. Let us denote
by $P_1$ and $P_2$ the two singular points
of $\Rr X_-$ over the two singular points of $\Cc A$.  Since
\[\begin{array}{ll}c_1(\Cc X_-)&=p^*\big(c_1(\Cc P(2,1,1))+\frac{1}{2}[\Cc A]\big)\\
&= \big(2+\frac{m}{2}\big) p^*(S_{\infty})+\big(2(k+1)+mk\big)p^*(F),\end{array}\]
$\Rr\widetilde{X}_-$ is orientable if and only if $4|m$.

The complex surface $\overline{P}=\Cc P(2k,1,1)/\conj$ is a
$\Qq$-homology sphere. Since the {\it Arnold surface}
$\mathcal{A}=\overline{A}\cup_{\partial} c(\Delta)$ is a 2-sphere,
it follows that $\overline{X}_-=\Cc X_-/\conj$ is a $\Qq$-homology 4-sphere.

If $4|m$, the surface $\Rr X_-$ is an oriented surface of genus 2 with two of its
nontrivial cycles pinched (see Figure \ref{genus2}); the two points $P_3$ and $P_4$
obtained in this way are Q-singularities with
\[
[\tb_{P_3}]=[\tb_{P_4}]=
\left[ \!\!\begin{array}{rr}
-\frac{2k}{4}&\frac{2k}{4}\\\frac{2k}{4}&-\frac{2k}{4}\end{array}\right]
\]
observed by Corollary \ref{contract2}, Section \ref{sec:6}.
\begin{figure}[h]
\begin{center}
$P_3$ \hspace{2.4cm} $P_4$

\vspace{-0.2cm}

$P_1$

\vspace{0.6cm}

$P_2$

\vspace{-1.75cm}

\epsfig{file=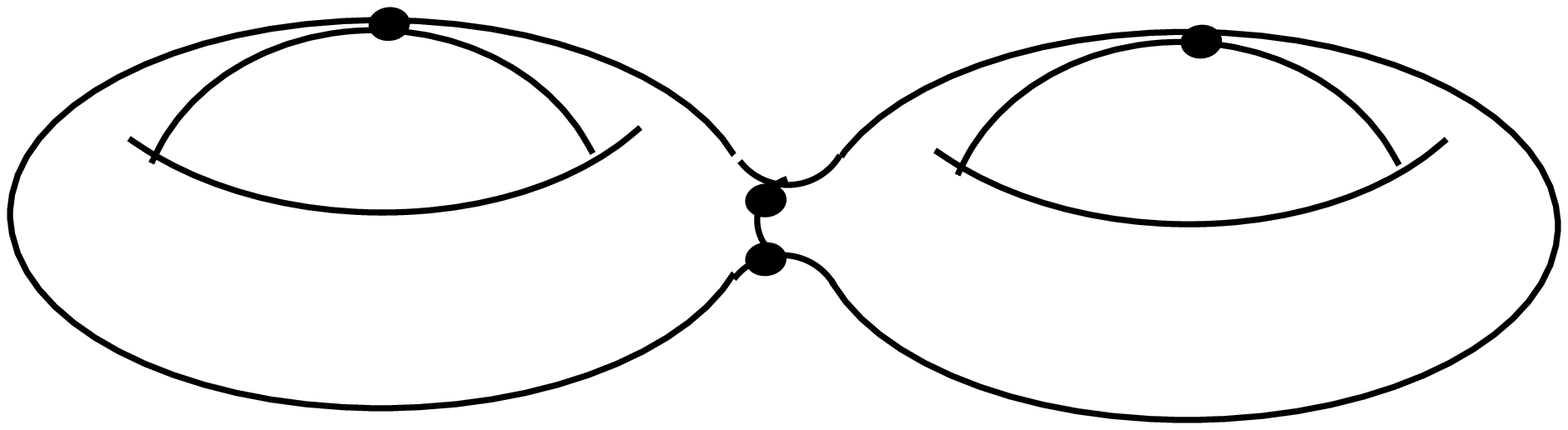,height=2cm,width=6cm} \caption{}
\label{genus2}

\end{center}
\end{figure}

\[
\begin{array}{ll}
\mbox{Hence, \hspace{1.5cm}  } 0&=\big\langle c'(\Rr X_-),c'(\Rr X_-)\big\rangle_{\overline{X}_-}\\
&=\sum_{i=1}^4\overline{\tb}_{P_i}(\Rr X_-,\Rr X_-)
-2\chi(\Rr X_--4\mbox{ points})\\
& =\overline{\tb}_{P_1} +\overline{\tb}_{P_2}
+0+0+-2\cdot(-2-2)\\
&=\overline{\tb}_{P_1} +\overline{\tb}_{P_2}+8.
\end{array}\] The proof for the first claim  now follows from Equation \ref{tbbar2rb}, Section \ref{sec:1}.

If $m\equiv 2 \mbox{ mod } 4$, then $\Rr\widetilde{X}_-$ is not orientable.
But as in the proof of Theorem \ref{calcul},
we observe that $\Rr X_-=\Rr\widetilde{X}_-/\widetilde{p}^{-1}(S_\infty)$ is orientable.
The two Q-singularities $P_3$ and $P_4$
obtained in this way have
\[
[\tb_{P_3}]=[\tb_{P_4}]=
\left[ \!\!\begin{array}{rr}
-\frac{2k}{4}&-\frac{2k}{4}\\-\frac{2k}{4}&-\frac{2k}{4}\end{array}\right]
\] so that
\[
\begin{array}{ll}
0 &=\sum_{i=1}^4\overline{\tb}_{P_i}(\Rr X_-,\Rr X_-)
-2\chi(\Rr X_--4\mbox{ points})\\
& =\overline{\tb}_{P_1} +\overline{\tb}_{P_2}
+2\big(-\frac{k}{2}-\frac{k}{2}-\frac{k}{2}-\frac{k}{2}\big)+-2\cdot(-2-2)\\
&=\overline{\tb}_{P_1} +\overline{\tb}_{P_2}-8k+8
\end{array}\]
completing the proof for the second claim.
\end{proof}

%
%

\end{document}